\documentclass[12pt]{article}
\usepackage{subeqn,amsfonts}
\newcommand{\be}{\begin{equation}}
\newcommand{\ee}{\end{equation}}
\newcommand{\bea}{\begin{eqnarray}}
\newcommand{\eea}{\end{eqnarray}}
\newcommand{\beo}{\begin{eqnarray*}}
\newcommand{\eeo}{\end{eqnarray*}}

\newcommand{\eps}{{\varepsilon}}
\newcommand{\cA}{{\cal A}}
\newcommand{\cG}{{\cal G}}
\newcommand{\cU}{{\cal U}}
\newcommand{\NN}{{\mathbb N}}
\newcommand{\ZZ}{{\mathbb Z}}

\newcommand{\CC}{{\mathbb C}}
\newcommand{\tad}[2]{{\Big(\raisebox{-2mm}{\shortstack{#1\\#2}}\Big)}}
\newcommand{\tadb}{\bigg(\tad{$a_1$}{$d_1$},\tad{$a_2$}{$d_2$}\bigg)}
\newcommand{\tadp}{\bigg(\tad{$1$}{$1$},\tad{$1$}{$1$}\bigg)}
\newcommand{\tadi}{\bigg(q^{-\eps_1}\tad{\phantom{-}$1$}{-$1$},
   q^{-\eps_2}\tad{\phantom{-}$1$}{-$1$}\bigg)}
\newcommand{\shalf}{\textstyle{\frac{1}{2}}}
\newcommand{\squart}{\textstyle{\frac{1}{4}}}

\newtheorem{thm}{Theorem}
\newtheorem{lem}{Lemma}

\begin{document}
\newpage
\pagestyle{empty}
\setcounter{page}{0}

\vfill
\vfill

\begin{center}

{\Large {\bf Classification of the quantum deformations

\vspace{5mm}

of the superalgebra $gl(1|1)$}}

\vspace{7mm}

{\large L. Frappat \footnote{On leave of absence from Laboratoire de 
Physique Th\'eorique ENSLAPP, Annecy, France.}, V. Hussin}

\vspace{4mm}

{\em Centre de Recherches Math\'ematiques, Universit\'e de Montr\'eal, Canada}

\vspace{7mm}

{\large G. Rideau}

\vspace{4mm}

{\em Laboratoire de Physique Th\'eorique et Math\'ematique, Universit\'e
Paris VII, France}

\end{center}

\vfill

\begin{abstract}
We present a classification of the possible quantum deformations of the
supergroup $GL(1|1)$ and its Lie superalgebra $gl(1|1)$. In each case, 
the (super)commutation relations and the Hopf structures are explicitly
computed. For each $R$ matrix, one finds two inequivalent coproducts
whether one chooses an unbraided or a braided framework while the 
corresponding structures are isomorphic as algebras. In the braided case, 
one recovers the classical algebra $gl(1|1)$ for suitable limits of the 
deformation parameters but this is no longer true in the unbraided case.
\end{abstract}

\def\abstractname{R\'esum\'e}
\begin{abstract}
Nous pr\'esentons une classification des d\'eformations quantiques du 
supergroupe $GL(1|1)$ et de sa superalg\`ebre $gl(1|1)$. Dans chaque cas,
les relations de (super)commutation et les structures de Hopf sont calcul\'ees
explicitement. Pour chaque matrice $R$, on trouve deux coproduits 
in\'equivalents selon que l'on choisit un sch\'ema tress\'e ou non, alors 
que les structures correspondantes sont isomorphes en tant qu'alg\`ebres. 
Dans le cas tress\'e, on retrouve l'alg\`ebre classique $gl(1|1)$ pour des
limites convenables des param\`etres de d\'eformation, mais ceci n'est
plus vrai dans le cas non tress\'e.
\end{abstract}

\vfill
\vfill

\rightline{CRM-2476}
\rightline{q-alg/9705024}
\rightline{May 1997}

\newpage
\pagestyle{plain}

\section{Introduction}

The method of $R$-matrix [1--5] for constructing 
quantum groups has already been generalized to quantum supergroups. For 
example, three non-equivalent such quantum supergroups have been derived 
recently for the fermionic oscillator group \cite{HLR96}. Another example 
deals with $GL(1|1)$, another four dimensional supergroup. The standard 
one-parameter deformation $GL_q(1|1)$ is well-known \cite{SVM90,CFFS90,LS92} 
and has been generalized to two parameters \cite{DW91,BBCH92}. An alternative 
deformation has also been derived \cite{Hla92}.

These last deformations are based on the choice of an $4 \times 4$ $R$-matrix 
which satisfies the constant quantum Yang-Baxter equation (YBE). A complete 
set of such solutions has been constructed \cite{Hie92} and may be the 
starting point for considering all possible continuous deformations of the 
linear group $GL(2)$ and  supergroup $GL(1|1)$. The extra conditions that 
have to be satisfied to this aim leads us to pick only solutions which are 
nonsingular $R$-matrices continuously related to some diagonal matrices.

It is already known that all the possible deformations of $GL(2)$ 
that possess a central determinant are given by the standard one 
\cite{Dri86,FRT90} and the non-standard (or ``Jordanian'') one 
\cite{Zak91,Ohn92}. Let us mention that they are both one-parameter 
deformations. Once the condition of central determinant is relaxed, we can 
show \cite{HLV97} that this ``Jordanian'' matrix contains two parameters 
and the computation of the quantum algebra dual to the quantum group is 
much more difficult and not known.

The quantum deformations of the group $GL(1|1)$ has until now not led 
to an exhaustive study. So the question addressed here is to give such 
a study and construct deformations both of the supergroup and 
superalgebra structures. Let us notice that our approach deals with 
deformations of superstructures with even parameters in comparison with
other recent approaches \cite{DPa96}.

A point which is important and has already been mentionned \cite{Maj95} is 
the fact that what distinguish the group and supergroup deformations is
that the corresponding R-matrices are continuous deformations
of the identity matrix in the first case and of the superidentity matrix 
(i.e. $\mbox{diag}(1,1,1,-1)$) in the second. 

While the paper will be concerned by the supergroup deformations, it is 
useful to present the necessary definitions for the usual group $GL(2)$ and 
point out the differences for the supergroup $GL(1|1)$.

So, let us consider the Lie group $G=GL(2)$, its Lie algebra $\cG=gl(2)$ 
with generators $A,B,C,D$ such that
\be
\Big[ A,B \Big] = - \Big[ D,B \Big] = B \,, \quad
\Big[ A,C \Big] = - \Big[ D,C \Big] = -C \,, \quad
\Big[ B,C \Big] = A-D \,, \quad \Big[ A,D \Big] = 0 \,, \label{defalg1}
\ee
and $\cU$ the universal enveloping algebra of $\cG$. The algebra 
$\cA = Fun(GL(2))$ is the associative unital algebra with generators 
$a,b,c,d$ that commute:
\be
[a,b] = [a,c] = [a,d] = [b,c] = [b,d] = [c,d] = 0 \,. \label{defgrp1}
\ee
The two algebras $\cU$ and $\cA$ can be endowed with a Hopf structure, each
element of $\cG\subset\cU$ being primitive for the comultiplication $\Delta$ 
(i.e. $\forall \ X \in \cG \,, \ \Delta(X) = X \otimes 1 + 1 \otimes X$) 
and the comultiplication $\Delta$ for $\cA$ is implied by the usual matrix
multiplication law: 
$\Delta T = T \, \raisebox{-4pt}{\shortstack{$\otimes$\\,}} \, T$ if
$T = \left( \begin{array}{cc} a&b \\ c&d \end{array} \right)$, that is
\be
\begin{array}{ll}
\Delta a = a \otimes a + b \otimes c \,, & 
\Delta b = a \otimes b + b \otimes d \,, \cr
\Delta c = c \otimes a + d \otimes c \,, & 
\Delta d = c \otimes b + d \otimes d \,. \cr
\end{array}
\label{coprod}
\ee
Moreover, there exists a nondegenerate bilinear form $\langle \ ,\ \rangle$
on $\cU \times \cA$ such that
\be
\begin{array}{ll}
\langle A,a^kd^lb^mc^n \rangle = k \delta_{m0} \delta_{n0} \,, &
\langle B,a^kd^lb^mc^n \rangle = \delta_{m1} \delta_{n0} \,, \cr
\langle C,a^kd^lb^mc^n \rangle = \delta_{m0} \delta_{n1} \,, &
\langle D,a^kd^lb^mc^n \rangle = l \delta_{m0} \delta_{n0} \,, \cr
\end{array}
\label{dual}
\ee
where $a^kd^lb^mc^n$ is any element of a Poincar\'e--Birkhoff--Witt 
basis of $\cA$ ($k,l,m,n \in \NN$).
Finally, the pairing $\langle \ ,\ \rangle$ satisfies
\be
\langle P_1 P_2,x \rangle = \langle m(P_1 \otimes P_2),x \rangle = 
\langle P_1 \otimes P_2,\Delta(x) \rangle \label{dualb}
\ee
and
\be
\langle \Delta(P),x \otimes y \rangle = \langle P,m(x \otimes y) \rangle = 
\langle P,xy \rangle \label{dualt}
\ee
where $P_1,P_2 \in \cU$, $x,y \in \cA$ and $m$ denotes the multiplication.
\\
The relations (\ref{dualb}) and (\ref{dualt}) make the Hopf algebras
$\cU$ and $\cA$ dual to each other.

\medskip

These definitions may be extended to the Lie supergroup $GL(1|1)$ and the Lie
superalgebra $gl(1|1)$. If we call as before its generators by $A,B,C,D$, then
we have
\bea
&& \Big[ A,B \Big] = - \Big[ D,B \Big] = B \,, \quad
\Big[ A,C \Big] = - \Big[ D,C \Big] = -C \,, \nonumber \\
&& \Big\{ B,C \Big\} = A+D \,, \quad \Big[ A,D \Big] = 
\Big\{ B,B \Big\} = \Big\{ C,C \Big\} = 0 \,. \label{defalg2}
\eea
Now the algebra $\cA = Fun(GL(1|1))$ with generators $a,b,c,d$ satisfies:
\be
[a,b] = [a,c] = [a,d] = [b,c] = 0 \,, \quad 
\{b,d\} = \{c,d\} = b^2 = c^2 = 0 \,. \label{defgrp2}
\ee

Deformations of the defining relations (\ref{defalg1})--(\ref{defgrp1}) 
or (\ref{defalg2})--(\ref{defgrp2}) are
provided by the Faddeev--Reshetikhin--Takhtajan formalism \cite{FRT90}. 
Let us define $T_1 = T \otimes I$, $T_2 = I \otimes T$. Then the deformations 
are given by
\be
R T_1 T_2 = T_2 T_1 R \,, \label{frt}
\ee
where $R$ is a $4 \times 4$ matrix that satisfies the quantum Yang-Baxter
equation (YBE)
\be
R_{12} R_{13} R_{23} = R_{23} R_{13} R_{12} \,, \label{ybe}
\ee
the last equation standing in $\mbox{End}(\CC^2) \otimes \mbox{End}(\CC^2) 
\otimes \mbox{End}(\CC^2)$.

As we said before, the $4 \times 4$ constant $R$-matrices satisfying the 
YBE have been classified in \cite{Hie92} and among them, the subset of non 
singular $R$-matrices splits into two different classes:
\\
$i$) the ones continuously connected to the identity matrix
$\mbox{diag}(1,1,1,1)$, which yield to quantum deformations of the group 
$GL(2)$: eq. (\ref{frt}) deforms the relations (\ref{defgrp1});
\\
$ii$) the ones continuously connected to the diagonal matrix 
$\mbox{diag}(1,1,1,-1)$, which yield to quantum deformations of the 
supergroup $GL(1|1)$: eq. (\ref{frt}) deforms the relations (\ref{defgrp2}).

In the first class, there is only two distinct deformations (one case has 
been discussed by Fronsdal et al. \cite{Fro93} and some work 
\cite{Zak91,Ohn92} has been done on the second case specializing in the 
one-parameter deformation).

Let us in the following concentrate on the second class of deformations, 
namely the ones of $gl(1|1)$.

\section{Deformations of the supergroup $GL(1|1)$}

The class of $R$-matrices satisfying the YBE and continuously connected to
$\mbox{diag}(1,1,1,-1)$ consists of three inequivalent matrices:
\bea
&& R_{2,2} = \left( \begin{array}{cccc} r & 0 & 0 & 0 \cr 
0 & 1 & r(1-q^{-1}) & 0 \cr 0 & 0 & r^2q^{-1} & 0 \cr 0 & 0 & 0 & -rq^{-1} \cr
\end{array} \right) \,, \\
&& \nonumber \\
&& \nonumber \\
&& R_{1,2} = \left( \begin{array}{cccc} 1 & 0 & 0 & r \cr 
0 & 1 & 1-q^{-1} & 0 \cr 0 & 0 & q^{-1} & 0 \cr 0 & 0 & 0 & -q^{-1} \cr
\end{array} \right) \,, \\
&& \nonumber \\
&& \nonumber \\
&& R_{1,1} = \frac{1}{r} \left( \begin{array}{cccc} 
s+1 & 0 & 0 & s \cr 0 & r & s & 0 \cr 0 & s & r & 0 \cr s & 0 & 0 & s-1 \cr
\end{array} \right) \,.
\eea
The first two matrices are really two-parameter matrices while the last one 
is a one-parameter matrix, the numbers $r,s$ being subject to the condition 
$r^2-s^2=1$ for the matrix $R_{1,1}$.

\medskip

The first case is already known, but in order to get a complete 
classification, we remind here the results.
The multiplication law between the generators of $\cA_{2,2}$ is given by
\be
\begin{array}{ll}
\bigg. ba - rab = 0 \,, & \quad rca - qac = 0 \,, \\
\bigg. bd + rdb = 0 \,, & \quad rcd + qdc = 0 \,, \\
\bigg. ad - da + r^{-1}(q-1) bc = 0 \,, & \quad r^2cb - qbc = 0 \,, \\
\bigg. b^2 = c^2 = 0 \,. \\
\end{array}
\ee
\setcounter{thm}{-1}
\begin{thm}\label{thmzero} {\rm \cite{BBCH92}}
The supercommutation relations for the dual algebra $\cU_{2,2}$, quantum 
deformation of the Lie superalgebra $gl(1|1)$ associated to the $R$-matrix 
$R_{2,2}$, are given by:
\beo
&& \Big[ A,D \Big] = 0 \,, \hspace{30mm} \Big\{ C,C \Big\} 
= \Big\{ B,B \Big\} = 0 \,, \\
&& \Big[ A,B \Big] = - \Big[ D,B \Big] = B \,, \hspace{9mm}
\Big[ A,C \Big] = - \Big[ D,C \Big] = - C \,, \\
&& \Big\{ B,C \Big\} = \frac{q^{A+D}-1}{q-1} \,,
\eeo
and the comultiplication structure by:
\beo
&& \Delta(A) = 1 \otimes A + A \otimes 1 \,, \qquad
\Delta(B) = 1 \otimes B + B \otimes (-1)^D r^{A+D} \,, \\
&& \Delta(D) = 1 \otimes D + D \otimes 1 \,, \qquad
\Delta(C) = 1 \otimes C + C \otimes (-1)^D 
\left(\frac{q}{r}\right)^{A+D} \,.
\eeo
\end{thm}

\section{The case $R_{1,2}$}

The multiplication law between the generators of $\cA_{1,2}$ is obtained from 
(\ref{frt}) with $R = R_{1,2}$ as:
\be
\begin{array}{ll}
\bigg. ba - ab + rq dc = 0 \,, & \quad ca - q ac = 0 \,, \\
\bigg. bd + db - rq ac = 0 \,, & \quad cd + q dc = 0 \,, \\
\bigg. ad - da - (1-q) bc = 0 \,, & \quad cb - q bc = 0 \,, \\
\bigg. (1+q)b^2 - rq(a^2-d^2) = 0 \,, & \quad c^2 = 0 \,. \\
\end{array} \label{multone}
\ee
The structure relations of the corresponding dual algebra $\cU_{1,2}$ will 
be obtained by computing the action of the (anti)commutators between the 
generators $A,B,C,D$ of $\cU_{1,2}$ on a Poincar\'e--Birkhoff--Witt basis 
of $\cA_{1,2}$. Such a basis is generated by the generic elements of the type 
$a^kd^lb^mc^n$ where $k,l \in \NN$ and $m,n \in \{0,1\}$ thanks to the two 
last relations of (\ref{multone}). Moreover, eq. (\ref{dualb}) requires the 
knowledge of the comultiplication $\Delta(a^kd^lb^mc^n)$. In the case under 
consideration, such a computation can be done directly and we have the 
following lemma:
\begin{lem}
\beo
&& \Delta(a^k) = a^k \otimes a^k + \frac{q^k-1}{q-1} a^{k-1}b \otimes a^{k-1}c
- rq^2 \frac{(q^k-1)(q^{k-1}-1)}{(q^2-1)(q-1)} a^{k-2}dc \otimes a^{k-1}c
\,, \\
&& \Delta(d^l) = d^l \otimes d^l + \frac{q^l-1}{q-1} d^{l-1}c \otimes d^{l-1}b
- rq^2 \frac{(q^l-1)(q^{l-1}-1)}{(q^2-1)(q-1)} d^{l-1}c \otimes ad^{l-2}c \,.
\eeo
\end{lem}
{\bf Proof:}
These relations are easily proved by recurrence on $k$ and $l$, using 
eqs. (\ref{coprod}) and (\ref{multone}).
\hfill \rule{5pt}{5pt}

\medskip

\noindent
One has, from the fact that $\Delta$ is an algebra homomorphism,
\bea
\Delta(a^kd^l) &=& a^kd^l \otimes a^kd^l 
+ q^l\frac{(q^k-1)(q^l-1)}{(q-1)^2} (a^{k-1}d^{l-1}bc \otimes 
a^{k-1}d^{l-1}bc  - rq a^kd^{l-1}c \otimes a^{k-1}d^lc) \nonumber \\
&& + \frac{q^l-1}{q-1} \left( a^kd^{l-1}c \otimes a^kd^{l-1}b - rq^2 
\frac{q^{l-1}-1}{q^2-1} a^kd^{l-1}c \otimes a^{k+1}d^{l-2}c \right) 
\nonumber \\
&& + q^l\frac{q^k-1}{q-1} \left( a^{k-1}d^lb \otimes a^{k-1}d^lc - rq^{l+2} 
\frac{q^{k-1}-1}{q^2-1} a^{k-2}d^{l+1}c \otimes a^{k-1}d^lc \right) \,.
\label{copad}
\eea
In the same way, one can deduce
\bea
\Delta(a^kd^lc) &=& a^kd^lc \otimes a^{k+1}d^l + a^kd^{l+1} \otimes a^kd^lc 
- \frac{q^l-1}{q-1}a^kd^lc \otimes a^kd^{l-1}bc \nonumber \\
&& + q^{l+1}\frac{q^k-1}{q-1}a^{k-1}d^lbc \otimes a^kd^lc \,, \label{copadc}
\eea
and finally
\bea
\Delta(a^kd^lb) &=& a^{k+1}d^l \otimes a^kd^lb + a^kd^lb \otimes a^kd^{l+1}
- \frac{q^l-1}{q-1} a^kd^{l-1}bc \otimes a^kd^lb
+ q^{l+1} \frac{q^k-1}{q-1} a^kd^lb \otimes a^{k-1}d^lbc \nonumber \\
&& \hspace{-15mm}
- rq^{l+2} \frac{q^k-1}{q^2-1} a^{k-1}(a^2-d^2)d^l \otimes a^{k-1}d^{l+1}c 
+ rq^2 \frac{q^l-1}{q^2-1} a^{k+1}d^{l-1}c \otimes a^k(a^2-d^2)d^{l-1} 
\nonumber \\
&& \hspace{-15mm} - rq^2 \frac{(q^l-1)(q^{l-1}-1)}{(q^2-1)(q-1)} 
a^{k+1}d^{l-1}c \otimes a^{k+1}d^{l-2}bc + rq^{2l+4} \frac{(q^k-1)(q^{k-1}-1)}
{(q^2-1)(q-1)} a^{k-2}d^{l+1}bc \otimes a^{k-1}d^{l+1}c \nonumber \\
&& \hspace{-15mm} + rq^2 \frac{q^{k+l+1}-1}{(q^2-1)(q-1)} 
\Big( (q^l-1) a^kd^{l-1}bc \otimes a^{k+1}d^{l-1}c 
- q^l (q^k-1) a^{k-1}d^{l+1}c \otimes a^{k-1}d^lbc \Big) \nonumber \\
&& \hspace{-15mm} + rq^{l+2} \frac{(q^k-1)(q^l-1)}{(q^2-1)(q-1)} 
\Big( a^kd^{l-1}bc \otimes a^{k-1}d^{l+1}c - q\, a^{k+1}d^{l-1}c \otimes 
a^{k-1}d^lbc \Big) \,, \label{copadb}
\eea
\bea
\Delta(a^kd^lbc) &=& a^{k+1}d^lc \otimes a^{k+1}d^lb 
- a^kd^{l+1}b \otimes a^kd^{l+1}c + a^kd^lbc \otimes a^{k+1}d^{l+1} 
+ a^{k+1}d^{l+1} \otimes a^kd^lbc \nonumber \\
&& \hspace{-15mm} 
+ \frac{q^{l+1}(q^{k+1}-1)-(q^{l+1}-1)}{q-1} a^kd^lbc \otimes a^kd^lbc
+ \frac{rq^2}{q^2-1} (q^l-1 - q^{l+1}(q^k-1)) a^{k+1}d^lc \otimes a^kd^{l+1}c
\nonumber \\
&& \hspace{-15mm} 
- rq^2 \frac{q^l-1}{q^2-1} a^{k+1}d^lc \otimes a^{k+2}d^{l-1}c  
+ rq^2 q^{l+1} \frac{q^k-1}{q^2-1} a^{k-1}d^{l+2}c \otimes a^kd^{l+1}c \,.
\label{copadbc}
\eea
Now we can enounce the following result:
\begin{thm}\label{thmalgone}
The supercommutation relations for the dual algebra $\cU_{1,2}$, quantum 
deformation of the Lie superalgebra $gl(1|1)$ associated to the $R$-matrix 
$R_{1,2}$, are given by (we have set $K = q^{A+D}$):
\beo
&& \Big[ A,D \Big] = 0 \,, \hspace{30mm} \Big\{ B,C \Big\} = 
\frac{q^{A+D}-1}{q-1} \,, \\
&& \Big[ A,B \Big] = - \Big[ D,B \Big] = B \,, \\
&& \Big[ A,C \Big] = - \Big[ D,C \Big] = - C - \frac{2rq}{q^2-1} \, 
(K-q)B \,, \\
&& \Big\{ C,C \Big\} = \frac{-2rq}{(q^2-1)(q-1)} \, (K-1)(K-q) 
\,, \qquad \Big\{ B,B \Big\} = 0 \,.
\eeo
\end{thm}
Note that the element $K$ is central in $\cU_{1,2}$: 
$\Big[ K,X \Big] = 0$ for $X \in \{A,B,C,D\}$.

\medskip

\noindent
{\bf Proof:} 
Using the formulae (\ref{copad}), (\ref{copadc}), (\ref{copadb}) and 
(\ref{copadbc}), we see that the non vanishing pairings are the following:
\bea
&& \langle BC+CB,a^kd^l \rangle = \frac{q^{k+l}-1}{q-1} \,, \nonumber \\
&& \langle AB-BA,a^kd^lb \rangle = -\langle DB-BD,a^kd^lb \rangle = 
1 \,, \nonumber \\
&& \langle AC-CA,a^kd^lc \rangle = -\langle DC-CD,a^kd^lc \rangle = 
-1 \,, \label{pair} \\
&& \langle AC-CA,a^kd^lb \rangle = -\langle DC-CD,a^kd^lb \rangle = 
-2rq^2 \frac{q^{k+l}-1}{q^2-1} \,, \nonumber \\
&& \langle C^2,a^kd^l \rangle = -rq^2 \frac{(q^{k+l}-1)(q^{k+l-1}-1)}
{(q^2-1)(q-1)} \,. \nonumber
\eea
To go from formulae (\ref{pair}) to the equations of Theorem \ref{thmalgone}, 
we need the following expressions:
\be
\begin{array}{l}
\langle A^n,a^kd^l \rangle = \langle \otimes_n A,
\Delta^{(n)}(a^kd^l) \rangle = k^n \\
\langle D^n,a^kd^l \rangle = \langle \otimes_n D,
\Delta^{(n)}(a^kd^l) \rangle = l^n
\end{array}
\ee
obtained from the coproduct (\ref{coprod}) and the multiplication law 
(\ref{multone}). It follows immediately that
\be
\langle q^A,a^kd^l \rangle = q^k \quad\mbox{and}\quad
\langle q^D,a^kd^l \rangle = q^l \,. \label{dualexp}
\ee
Moreover, one has from eq. (\ref{copadb}) (note the shift in the exponential !)
\be
\langle q^{A+D-1}B,a^kd^lb \rangle = q^{k+l} \,. \label{eq23}
\ee
Then comparing eqs. (\ref{pair}), (\ref{dualexp}) and (\ref{eq23}), we get 
the commutation relations of Theorem \ref{thmalgone}.
\hfill \rule{5pt}{5pt}

\medskip

We want now to determine the comultiplication structure on $\cU_{1,2}$.
The duality relation (\ref{dualt}) applied on the generic elements 
$a^kd^lb^mc^n$ and $a^{k'} d^{l'} b^{m'} c^{n'}$ of the 
Poincar\'e--Birkhoff--Witt basis of $\cA_{1,2}$ reads as
\be
\langle \Delta(P),a^kd^lb^mc^n \otimes a^{k'}d^{l'}b^{m'}c^{n'} \rangle
= \langle P,m(a^kd^lb^mc^n \otimes a^{k'}d^{l'}b^{m'}c^{n'}) \rangle 
= \langle P,a^kd^lb^mc^n a^{k'}d^{l'}b^{m'}c^{n'} \rangle \,.
\label{dualpwb}
\ee
If $\Delta(P) = P_{(1)} \otimes P_{(2)}$ in Sweedler's notation, one has
\be
\langle \Delta(P),a^kd^lb^mc^n \otimes a^{k'}d^{l'}b^{m'}c^{n'} \rangle
= \langle P_{(1)},a^kd^lb^mc^n \rangle 
\langle P_{(2)},a^{k'}d^{l'}b^{m'}c^{n'} \rangle \,. \label{sweed}
\ee
{From} the knowledge of $\langle P,a^kd^lb^mc^na^{k'}d^{l'}b^{m'}c^{n'} 
\rangle$ as a function of $k,l,m,n,k',l',m',n'$, and the duality
relations (\ref{dual}), one can then deduce the possible $P_{(1)}$ and 
$P_{(2)}$ for any generator $P$ of the dual algebra.

{From} formula (\ref{dualpwb}), one has to compute the action of any generator 
of the algebra $\cU_{1,2}$ on a generic element $a^kd^lb^mc^na^{k'}d^{l'} 
b^{m'}c^{n'}$ where $m,n,m',n'=0$ or $1$. Using the multiplication law 
(\ref{multone}), it is possible to reorder this generic element with respect 
to the ordering $adbc$ given by the duality relations (\ref{dual}).
To this aim, we need the following lemma (reordering formulae):
\begin{lem}\label{lemtwo}
\beo
a^kd^la^{k'}d^{l'} &=& a^{k+k'}d^{l+l'} + q^{l'}\frac{(q^{k'}-1)(q^l-1)}
{q-1} a^{k+k'-1}d^{l+l'-1}bc \,, \\
a^kd^la^{k'}d^{l'}b &=& a^{k+k'}d^{l+l'}b + \frac{rq^2}{q^2-1}
q^{l'}(q^{k'}-1)(q^l-1) a^{k+k'-1}(a^2-d^2)d^{l+l'-1}c \,, \\
a^kd^la^{k'}d^{l'}c &=& a^{k+k'}d^{l+l'}c \,, \\
a^kd^lba^{k'}d^{l'} &=& (-1)^{l'} a^{k+k'}d^{l+l'}b 
+ (-1)^{l'} \frac{rq^2}{q^2-1} q^{l'}(q^{k'}-1)(q^l-1) 
a^{k+k'-1}(a^2-d^2)d^{l+l'-1}c \\
&& - (-1)^{l'} \frac{rq}{q-1} (q^{l'}-1) a^{k+k'+1}d^{l+l'-1}c  
- (-1)^{l'} \frac{rq^{l'+1}}{q-1} (q^{k'}-1) a^{k+k'-1}d^{l+l'+1}c \,, \\
a^kd^lca^{k'}d^{l'} &=& (-1)^{l'} q^{k'+l'} a^{k+k'}d^{l+l'}c \,, \\
a^kd^lba^{k'}d^{l'}b &=& (-1)^{l'} \frac{rq}{q+1} a^{k+k'}(a^2-d^2)d^{l+l'} 
+ (-1)^{l'} rq(q^{l+l'}-1) a^{k+k'+1}d^{l+l'-1}bc \\
&& - (-1)^{l'} \frac{rq^2}{q-1}(q^{l'}-1) a^{k+k'+1}d^{l+l'-1}bc
- (-1)^{l'} \frac{rq^2}{q-1}q^{l'}(q^{k'}-1) a^{k+k'-1}d^{l+l'+1}bc \\
&& + (-1)^{l'} \frac{rq^3}{q^2-1} q^{l'}(q^{k'}-1)(q^l-1)
a^{k+k'-1}(a^2-d^2)d^{l+l'-1}bc \,, \\
a^kd^lba^{k'}d^{l'}c &=& (-1)^{l'} a^{k+k'}d^{l+l'}bc \,, \\
a^kd^lca^{k'}d^{l'}b &=& (-1)^{l'} q^{k'+l'+1} a^{k+k'}d^{l+l'}bc \,, \\
a^kd^lbca^{k'}d^{l'} &=& q^{k'+l'} a^{k+k'}d^{l+l'}bc \,, \\
a^kd^la^{k'}d^{l'}bc &=& a^{k+k'}d^{l+l'}bc \,, \\
a^kd^lba^{k'}d^{l'}bc &=& (-1)^{l'} \frac{rq}{q+1} 
a^{k+k'}(a^2-d^2)d^{l+l'}c \,, \\
a^kd^lbca^{k'}d^{l'}b &=& \frac{rq^2}{q+1} q^{k'+l'} 
a^{k+k'}(a^2-d^2)d^{l+l'}c \,, \\
a^kd^lca^{k'}d^{l'}c &=& a^kd^lca^{k'}d^{l'}bc \ =\ a^kd^lbca^{k'}d^{l'}c 
\ =\ a^kd^lbca^{k'}d^{l'}bc \ =\ 0 \,. 
\eeo
\end{lem}
{\bf Proof:}
The proof of the lemma is straightforward and is done by recurrence on 
$k,l,k',l'$ from eq. (\ref{multone}).
\hfill \rule{5pt}{5pt}
\begin{thm}\label{thmcopone}
The comultiplication $\Delta$ of the algebra $\cU_{1,2}$ is given by:
\beo
&& \Delta(A) = 1 \otimes A + A \otimes 1 + \frac{2rq}{q+1} B \otimes 
(-1)^D B \,, \\
&& \Delta(B) = 1 \otimes B + B \otimes (-1)^D \,, \\
&& \Delta(C) = 1 \otimes C + C \otimes (-1)^D K - \frac{rq}{q-1} B 
\otimes (-1)^D (K-1) \,, \\
&& \Delta(D) = 1 \otimes D + D \otimes 1 - \frac{2rq}{q+1} B \otimes 
(-1)^D B \,.
\eeo
\end{thm}
Let us remark that the first and last equations of Theorem \ref{thmcopone} 
imply that $\Delta(K) = K \otimes K$.

\medskip

\noindent
{\bf Proof:} 
It follows immediatly from Lemma \ref{lemtwo} that
\bea
&& \langle \Delta(A), a^kd^l \otimes a^{k'}d^{l'} \rangle = k+k' \,, \qquad
\ \langle \Delta(A), a^kd^lb \otimes a^{k'}d^{l'}b \rangle = (-1)^{l'} 
\frac{2rq}{q+1} \,, \nonumber \\
&& \langle \Delta(B), a^kd^l \otimes a^{k'}d^{l'}b \rangle = 1 \,, \qquad 
\qquad \langle \Delta(B), a^kd^lb \otimes a^{k'}d^{l'} \rangle = (-1)^{l'} \,, 
\nonumber \\
&& \langle \Delta(C), a^kd^l \otimes a^{k'}d^{l'}c \rangle = 1 \,, \qquad 
\qquad \langle \Delta(C), a^kd^lc \otimes a^{k'}d^{l'} \rangle = (-1)^{l'} 
q^{k'+l'} \,, \\
&& \langle \Delta(C), a^kd^lb \otimes a^{k'}d^{l'} \rangle = 
-(-1)^{l'} rq \frac{q^{k'+l'}-1}{q-1} \,, \nonumber \\
&& \langle \Delta(D), a^kd^l \otimes a^{k'}d^{l'} \rangle = l+l' \,, \qquad
\ \langle \Delta(D), a^kd^lb \otimes a^{k'}d^{l'}b \rangle = -(-1)^{l'} 
\frac{2rq}{q+1} \,, \nonumber
\eea
and all other possible terms vanish. These last relations then imply 
Theorem \ref{thmcopone} by using the duality relations (\ref{dual}), 
(\ref{dualt}) and (\ref{dualexp}). This achieves the proof.
\hfill \rule{5pt}{5pt}

\section{The case $R_{1,1}$}

The multiplication law between the generators of $\cA_{1,1}$ is given by the 
relation (\ref{frt}) with $R = R_{1,1}$. One obtains:
\be
\begin{array}{ll}
\bigg. ba - r ab + s dc = 0 \,, & \quad ca - r ac + s db = 0 \,, \\
\bigg. bd + r db - s ac = 0 \,, & \quad cd + r dc - s ab = 0 \,, \\
\bigg. ad-da = 0 \,, & \quad bc-cb = 0 \,, \\
\bigg. b^2 = c^2 = \shalf s (a^2-d^2) \,. \\
\end{array} \label{multtwo}
\ee
As before a Poincar\'e--Birkhoff--Witt basis of $\cA_{1,1}$ is given by 
$a^kd^lb^mc^n$ where $k,l \in \NN$ and $m,n \in \{0,1\}$ thanks to the last 
relation of (\ref{multtwo}). The computation of the coproduct 
$\Delta(a^kd^lb^mc^n)$ where $k,l \in \NN$ and $m,n \in \{0,1\}$ is much 
more involved than in the case $\cU_{1,2}$ because the 
multiplication law (\ref{multtwo}) does not allow to compute directly the 
quantities $\Delta(a^kd^lb^mc^n)$. Instead, one has to solve many recursion
formulae for $\Delta(a^kd^l)$ in order to produce the desired results (see
the Appendix).

\medskip

\begin{thm}\label{thmalgtwo}
The supercommutation relations for the dual algebra $\cU_{1,1}$, quantum 
deformation of the Lie superalgebra $gl(1|1)$ associated to the $R$-matrix 
$R_{1,1}$, are given by:
\beo
&&\Big[ A,D \Big] = 0 \,, \hspace{30mm}
\Big\{ B,C \Big\} = \frac{1}{2} \left( \frac{K^2-1}{q^2-1} + 
\frac{K^{-2}-1}{q^{-2}-1} \right) \,, \\
&&\Big[ A,B \Big] = - \Big[ D,B \Big] = \shalf B 
+ \squart (q^{-2}K^2+q^2K^{-2})B + \squart (q^{-2}K^2-q^2K^{-2})C \,, \\
&&\Big[ A,C \Big] = - \Big[ D,C \Big] = - \shalf C 
- \squart (q^{-2}K^2+q^2K^{-2})C - \squart (q^{-2}K^2-q^2K^{-2})B \,, \\
&&\Big\{ B,B \Big\} = \Big\{ C,C \Big\} = 
-\frac{1}{2} \left( \frac{K^2-1}{q^2-1} - 
\frac{K^{-2}-1}{q^{-2}-1} \right) \,.
\eeo
Since $r^2-s^2=1$, we have set for convenience $r = \shalf (q+q^{-1})$, 
$s = \shalf (q-q^{-1})$ and $K$ is defined by $K = q^{A+D}$.
\end{thm}
Note again that the element $K$ is central in $\cU_{1,1}$.

\medskip

\noindent
{\bf Proof:} 
The interested reader will find the details in the Appendix.
\hfill \rule{5pt}{5pt}

\medskip

\begin{thm}\label{thmcoptwo}
The comultiplication $\Delta$ of the algebra $\cU_{1,1}$ is given by:
\beo
&& \Delta(A) = 1 \otimes A + A \otimes 1 \nonumber \\
&& \hspace {15mm} + \squart (q-q^{-1}) \Big( (B-C) \otimes (-1)^D q^{-1} K 
(B+C) + (B+C) \otimes (-1)^D q K^{-1} (B-C) \Big) \,, \\
&& \Delta(B) = 1 \otimes B + \shalf (B-C) \otimes (-1)^D K
+ \shalf (B+C) \otimes (-1)^D K^{-1} \,, \\
&& \Delta(C) = 1 \otimes C - \shalf (B-C) \otimes (-1)^D K
+ \shalf (B+C) \otimes (-1)^D K^{-1} \,, \\
&& \Delta(D) = 1 \otimes D + D \otimes 1 \nonumber \\
&& \hspace {15mm} - \squart (q-q^{-1}) \Big( (B-C) \otimes (-1)^D q^{-1} K 
(B+C) + (B+C) \otimes (-1)^D q K^{-1} (B-C) \Big) \,. 
\eeo
\end{thm}
Again the first and last equations of Theorem \ref{thmcoptwo} imply that 
$\Delta(K) = K \otimes K$.

\medskip

\noindent
{\bf Proof:} The proof of Theorem \ref{thmcoptwo} stands along the same 
lines than the proof of Theorem \ref{thmcopone}.
{From} formula (\ref{dualpwb}), one computes the action of any generator 
of the algebra $\cU_{1,1}$ on a generic element $a^kd^lb^mc^na^{k'}d^{l'} 
b^{m'}c^{n'}$ where $m,n,m',n' \in \{0,1\}$. All that remains to do is
to reorder this generic element with respect to the ordering $adbc$ given 
by the duality relations (\ref{dual}). The reordering formulae are much 
simpler than in the $\cU_{1,2}$ case. Indeed, we have from eq. (\ref{multtwo})
\be
(b \pm c) a^kd^l = (ra \mp sd)^k (\pm sa - rd)^l (b \pm c) \,,
\ee
hence
\be
\begin{array}{l} 
\bigg. b a^kd^l = \shalf \Big((ra-sd)^k (sa-rd)^l (b+c) 
+ (ra+sd)^k (-sa-rd)^l (b-c) \Big) \equiv \xi^b_{kl} \,, \\
\bigg. c a^kd^l = \shalf \Big((ra-sd)^k (sa-rd)^l (b+c) 
- (ra+sd)^k (-sa-rd)^l (b-c) \Big) \equiv \xi^c_{kl} \,.
\end{array} \label{reord}
\ee
Then, for any $X \in \{A,B,C,D\}$, one has
\bea
\langle \Delta(X),a^kd^lb^mc^n \otimes a^{k'}d^{l'}b^{m'}c^{n'} \rangle 
&=& \langle X,a^kd^lb^mc^n a^{k'}d^{l'}b^{m'}c^{n'} \rangle 
\nonumber \\
&=& \left\{ \begin{array}{ll} 
\Big. \langle X,a^{k+k'}d^{l+l'}b^{m'}c^{n'} \rangle & \mbox{if } m=n=0 \\
\Big. \shalf \langle X,a^{k+k'}d^{l+l'} \xi^b_{k'l'} b^{m'}c^{n'} \rangle 
& \mbox{if } m=1,n=0 \\
\Big. \shalf \langle X,a^{k+k'}d^{l+l'} \xi^c_{k'l'} b^{m'}c^{n'} \rangle 
& \mbox{if } m=0,n=1 \\
\Big. \shalf \langle X,a^{k+k'}d^{l+l'} \xi^b_{k'l'} \xi^c_{k'l'} b^{m'} 
c^{n'} \rangle & \mbox{if } m=1,n=1 \\
\end{array} \right. \label{eqcop} 
\eea
It follows immediatly from eqs. (\ref{reord}) and (\ref{eqcop}) that
\begin{subequations}
\bea
&& \langle \Delta(A), a^kd^l \otimes a^{k'}d^{l'} \rangle = k+k' 
\,, \nonumber \\
&& \langle \Delta(A), a^kd^lb \otimes a^{k'}d^{l'}b \rangle = 
\langle \Delta(A), a^kd^lc \otimes a^{k'}d^{l'}c \rangle = 
\squart (q-q^{-1}) (-1)^{l'} (q^{k'+l'}+q^{-k'-l'}) 
\,, \\
&& \langle \Delta(A), a^kd^lb \otimes a^{k'}d^{l'}c \rangle = 
\langle \Delta(A), a^kd^lc \otimes a^{k'}d^{l'}b \rangle = 
-\squart (q-q^{-1}) (-1)^{l'} (q^{k'+l'}-q^{-k'-l'}) 
\,, \nonumber \\
&& \nonumber \\
&& \langle \Delta(B), a^kd^l \otimes a^{k'}d^{l'}b \rangle = 
\langle \Delta(C), a^kd^l \otimes a^{k'}d^{l'}c \rangle = 1 \,, \nonumber \\
&& \langle \Delta(B), a^kd^lb \otimes a^{k'}d^{l'} \rangle = 
\langle \Delta(C), a^kd^lc \otimes a^{k'}d^{l'} \rangle = 
\shalf (-1)^{l'} (q^{k'+l'}+q^{-k'-l'}) \,, \\
&& \langle \Delta(B), a^kd^lc \otimes a^{k'}d^{l'} \rangle = 
\langle \Delta(C), a^kd^lb \otimes a^{k'}d^{l'} \rangle = 
-\shalf (-1)^{l'} (q^{k'+l'}-q^{-k'-l'}) \,, \nonumber \\
&& \nonumber \\
&& \langle \Delta(D), a^kd^l \otimes a^{k'}d^{l'} \rangle = l+l' 
\,, \nonumber \\
&& \langle \Delta(D), a^kd^lb \otimes a^{k'}d^{l'}b \rangle = 
\langle \Delta(D), a^kd^lc \otimes a^{k'}d^{l'}c \rangle = 
-\squart (q-q^{-1}) (-1)^{l'} (q^{k'+l'}+q^{-k'-l'}) 
\,, \\
&& \langle \Delta(D), a^kd^lb \otimes a^{k'}d^{l'}c \rangle = 
\langle \Delta(D), a^kd^lc \otimes a^{k'}d^{l'}b \rangle = 
\squart (q-q^{-1}) (-1)^{l'} (q^{k'+l'}-q^{-k'-l'}) 
\,, \nonumber
\eea
\end{subequations}
and all other possible terms vanish. These last relations then imply 
Theorem \ref{thmcoptwo} by using the duality relations (\ref{dual}), 
(\ref{dualt}) and (\ref{prf23}), (\ref{prf26}) (for these relations, see 
the last $\clubsuit$ item of the Appendix). This achieves the proof.
\hfill \rule{5pt}{5pt}

\section{Braided structures}

In the case of the standard deformed superalgebra $gl(1|1)$, it is known
that there exist two different Hopf algebras $\cU_q[gl(1|1)]$ and
$\cU_q[gl(1|1)]'$, the two structures being isomorphic as algebras but
exhibiting two distinct Hopf structures. The former admits $gl(1|1)$ as
classical limit when $q \rightarrow 1$ while such a limit does not exist for 
the latter, $\cU_q[gl(1|1)]'$ being related to $\cU_i[sl(2,\CC)]$ at a root 
of unity ($i^2=-1$). A similar behaviour was proved in \cite{BBCH92} for the
$\cU_{2,2}$ case. This is a general feature as we will see below.

\medskip

The existence of two inequivalent Hopf structures is related to the fact that
one can choose a braided or an unbraided framework. 

In the unbraided case, the deformations of relations (\ref{defgrp2}) are 
given by (\ref{frt}): one finds the results stated in the previous sections. 
As can be seen from Theorems \ref{thmzero}, \ref{thmcopone} and 
\ref{thmcoptwo}, the corresponding deformations $\cU_{2,2}$, $\cU_{1,2}$, 
$\cU_{1,1}$ do not admit the classical superalgebra $gl(1|1)$ as a limit
for suitable values of the deformation parameters (it is clear from the
comultiplication formulae that $(-1)^D$ does not reduce to unity in such
a limit).

In the braided case, one has to introduce a ``braiding matrix'' chosen here
as the superidentity matrix $\mbox{diag}(1,1,1,-1)$. The braided version of
(\ref{frt}) reads as:
\be
R \hat{T}_1 \hat{T}_2 = \hat{T}_2 \hat{T}_1 R \,, \label{frtb}
\ee
where $\hat{T}_i = \eta T_i$ $(i=1,2)$.
\\
For $R=\mbox{diag}(1,1,1,-1)$, the generators of the algebra 
$\cA = Fun(GL(1|1))$ satisfy now (compare with relations (\ref{defgrp2}); 
note that the relations (\ref{defgrpb}) are consistent with a natural 
$\ZZ_2$-gradation with the assignment $a,d$ even and $b,c$ odd):
\be
[a,b] = [a,c] = [a,d] = [b,d] = [c,d] = 0 \,, \quad 
\{b,c\} = b^2 = c^2 = 0 \,. \label{defgrpb}
\ee

When the $R$-matrix is not trivial, it is easy to compute the modified 
multiplication laws for the cases $R_{2,2}$, $R_{1,2}$, $R_{1,1}$ 
corresponding to the deformations of (\ref{defgrpb}). One finds:
\begin{subequations}
\label{ab}
\bea
\mbox{for $\cA_{2,2}$:} \quad &&
\begin{array}{ll}
ba - rab = 0 \,, & \quad rca - qac = 0 \,, \\
bd - rdb = 0 \,, & \quad rcd - qdc = 0 \,, \\
ad - da + r(1-q^{-1})cb = 0 \,, & \quad r^2cb + qbc = 0 \,, \\
b^2 = c^2 = 0 \,. \\
\end{array} \\ && \nonumber \\
\mbox{for $\cA_{1,2}$:} \quad && 
\begin{array}{ll}
ba - ab + rq dc = 0 \,, & \quad ca - q ac = 0 \,, \\
bd - db + rq ac = 0 \,, & \quad cd - q dc = 0 \,, \\
ad - da + (1-q) bc = 0 \,, & \quad cb + q bc = 0 \,, \\
(1+q)b^2 - rq(a^2-d^2) = 0 \,, & \quad c^2 = 0 \,. \\
\end{array} \\ && \nonumber \\
\mbox{for $\cA_{1,1}$:} \quad &&
\begin{array}{ll}
ba - r ab + s dc = 0 \,, & \quad ca - r ac + s db = 0 \,, \\
bd - r db + s ac = 0 \,, & \quad cd - r dc + s ab = 0 \,, \\
ad - da = 0 \,, & \quad bc + cb = 0 \,, \\
b^2 = -c^2 = \shalf s (a^2-d^2) \,. \\
\end{array}
\eea
\end{subequations}
One can convince oneself, although it requires some work, that the 
(super)commutation relations of the corresponding dual algebras 
$\cU_{2,2}$, $\cU_{1,2}$, $\cU_{1,1}$ are unchanged. In this respect
the relations (\ref{ab}) just express the original algebras 
in a different basis. 
However, the Hopf structures are {\em not equivalent} to the ones presented 
in the previous sections. One finds the following results for the 
comultiplication (compare with Theorems \ref{thmzero}, \ref{thmcopone} and 
\ref{thmcoptwo}):
\begin{subequations}
\label{ub}
\bea
&\mbox{for $\cU_{2,2}$:}& \\
&&\Delta(A) = 1 \otimes A + A \otimes 1 \,, \qquad
\Delta(B) = 1 \otimes B + B \otimes r^{A+D} \,, \nonumber \\
&&\Delta(D) = 1 \otimes D + D \otimes 1 \,, \qquad
\Delta(C) = 1 \otimes C + C \otimes \left(\frac{q}{r}\right)^{A+D} 
\,. \nonumber \\
&& \nonumber \\
&\mbox{for $\cU_{1,2}$:}& \\
&&\Delta(A) = 1 \otimes A + A \otimes 1 + \frac{2rq}{q+1} B \otimes B \,,
\qquad \Delta(D) = 1 \otimes D + D \otimes 1 - \frac{2rq}{q+1} B \otimes B 
\,, \nonumber \\
&&\Delta(B) = 1 \otimes B + B \otimes 1 \,, \qquad
\Delta(C) = 1 \otimes C + C \otimes q^{A+D} - \frac{rq}{q-1} B \otimes 
(q^{A+D}-1) \,. \nonumber \\
&& \nonumber \\
&\mbox{for $\cU_{1,1}$:}& \\
&&\Delta(A) = 1 \otimes A + A \otimes 1 \nonumber \\
&& \hspace{15mm} + \squart (q-q^{-1}) \Big( (B-C) \otimes q^{A+D-1} (B+C) 
+ (B+C) \otimes q^{-A-D+1} (B-C) \Big) \,, \nonumber \\
&&\Delta(B) = 1 \otimes B + \shalf (B+C) \otimes q^{A+D} + \shalf (B-C) 
\otimes q^{-A-D} \,, \nonumber \\
&&\Delta(C) = 1 \otimes C + \shalf (B+C) \otimes q^{A+D} - \shalf (B-C) 
\otimes q^{-A-D} \,, \nonumber \\
&&\Delta(D) = 1 \otimes D + D \otimes 1 \nonumber \\
&& \hspace{15mm} - \squart (q-q^{-1}) \Big( (B-C) \otimes q^{A+D-1} (B+C) 
+ (B+C) \otimes q^{-A-D+1} (B-C) \Big) \,. \nonumber
\eea
\end{subequations}
Notice that the $q$-deformed superalgebras $\cU_{2,2}$, $\cU_{1,2}$ and
$\cU_{1,1}$ are now endowed with a super-Hopf structure, the comultiplication
$\Delta$ and the tensor product being $\ZZ_2$-graded, this last one satisfying
\be
(X_1 \otimes Y_1)(X_2 \otimes Y_2) = (-1)^{\deg Y_1.\deg X_2} 
(X_1 X_2 \otimes Y_1 Y_2) \,,
\ee
the $\ZZ_2$-gradation being defined by setting $\deg A= \deg D = 0$ and
$\deg B = \deg C = 1$.

\medskip

It is easy to see that the relations (\ref{ab}) and (\ref{ub}) lead
to the classical $GL(1|1)$ and $gl(1|1)$, endowing the superalgebra $gl(1|1)$
with a primitive comultiplication for suitable limits of the deformation
parameters: $r,q \rightarrow 1$ for the (2,2) case, $r \rightarrow 0, 
q \rightarrow 1$ for the (1,2) case and $q \rightarrow 1$ (or $r \rightarrow 1,
s \rightarrow 0)$ for the (1,1) case. Finally, the standard deformed
superalgebra $U_r[gl(1|1)]$ can be obtained by taking $q=r^2$ in the case 
$\cU_{2,2}$.

\section{Conclusion}

Starting with a two-dimensional representation of the supergroup $GL(1|1)$
we have been able to exhibit three types of continuous deformations of 
both the supergroup and superalgebra structures. These are based on the 
$R$-matrix method where $R$ satisfies the YBE. Two of the three types are
new with respect to preceding approches of the same question.

It is remarkable to notice that these results coincide, at the algebra
level, with those occuring in the fermionic oscillator quantum group approach
\cite{HLR96}. Indeed, the algebra corresponding to this fermionic oscillator 
appears to be isomorphic to $gl(1|1)$. We started with a three dimensional
representation of the corresponding group structure and obtained, with 
$9 \times 9$ $R$-matrices satisfying a weak version of YBE, three non 
isomorphic deformations of the superalgebra $gl(1|1)$ which can be compared 
with the ones obtained in this paper.

\medskip

For $\cU_{2,2}$ the correspondence is immediate and this superalgebra is 
relatied to the type III fermionic oscillator quantum superalgebra. 

For $\cU_{1,2}$, the change of basis $A'=A$, $D'=D$, $\displaystyle C' = C 
+ \Big(\frac{rq}{q^2-1}(K-q) + \frac{r(p-1)}{2p^2}(K-1) \Big)B$ and 
$\displaystyle B'=\frac{q-1}{p}B$ leads to (with $K=q^{A+D}$ and $p = \ln q$):
\beo
&& \Big[ A',D' \Big] = 0 \,, \hspace{30mm} \Big\{ B',C' \Big\} = 
\frac{K-1}{p} \,, \\
&& \Big[ A',B' \Big] = - \Big[ D',B' \Big] = B' \,, \\
&& \Big[ A',C' \Big] = - \Big[ D',C' \Big] = - C' - \frac{r}{p} 
\, (K-1)B' \,, \\
&& \Big\{ C',C' \Big\} = -\frac{r}{p^2} \, (K-1)^2
\,, \qquad \Big\{ B',B' \Big\} = 0 \,,
\eeo
which is related to the type II fermionic oscillator quantum superalgebra.

Finally, for $\cU_{1,1}$, we notice that a more natural basis would be
(with $K=q^{A+D}$):
\[
A' = q\frac{K^2+1}{K^2+q^2} \, A \,, \quad 
D' = q\frac{K^2+1}{K^2+q^2} \, D \,, \quad
B' = (1-q^{-2})^{1/2} \, (B+C) \,, \quad
C' = (q^2-1)^{1/2} \, (B-C) \,,
\]
since it gives
\beo
&& \Big[ A',D' \Big] = 0 \,, \hspace{30mm} \Big\{ B',C' \Big\} = 0 \,, \\
&& \Big[ A',B' \Big] = - \Big[ D',B' \Big] = \shalf (1+K^{-2}) C' \,, \\
&& \Big[ A',C' \Big] = - \Big[ D',C' \Big] = \shalf(1+K^2) B' \,, \\
&& \Big\{ B',B' \Big\} = 2(1-K^{-2}) \,, \qquad \Big\{ C',C' \Big\} = 
2(1-K^2) \,.
\eeo
This last structure is easily seen to be equivalent to the type I fermionic
oscillator quantum superalgebra (which clearly is a one-parameter deformation).

\subsubsection*{Acknowledgements}

The research of V.H. is partially supported by research grants from NSERC
of Canada and FCAR du gouvernement du Qu\'ebec.
L.F. is indebted to Centre de Recherches Math\'ematiques of Universit\'e
de Montr\'eal for its kind invitation and support.

\clearpage

\renewcommand{\thesection}{\Alph{section}}
\setcounter{section}{1}
\newtheorem{lemap}{Lemma}[section]
\section*{Appendix: Proof of theorem \ref{thmalgtwo}}

As stated above, the evaluation of the action of the generators of $\cU_{1,1}$
on the generic elements $a^kd^lb^mc^n$ of a Poincar\'e--Birkhoff--Witt basis 
of $\cA_{1,1}$ requires the calculation of the coproduct of such an element. 
Let us define 
\be
\Delta(a^kd^l) = \sum_{i,j,i',j' \in\{0,1\}} \Delta^{kl}_{ij,i'j'} 
\tadb \ b^ic^j \otimes b^{i'}c^{j'} \,. \label{prf1}
\ee
where the quantities $\Delta^{kl}_{ij,i'j'}\tadb$ are polynomials in the 
formal variables $a_1 = a \otimes 1$, $a_2 = 1 \otimes a$, 
$d_1 = d \otimes 1$, $d_2 = 1 \otimes d$.
\\
{From} the product formula (\ref{dualb}) and the duality relations 
(\ref{dual}), it is clear that the evaluation of the commutators between
the generators of $\cU_{1,1}$ on $a^kd^lb^mc^n$ is nothing but linear 
combinations of the polynomials $\Delta^{kl}_{ij,i'j'}\tadb$ and their 
derivatives for special values of the variables $a_1,a_2,d_1,d_2$. 
More precisely, if $P(a,d)$ is a polynomial in the variables $(a,d)$, the 
duality relations (\ref{dual}) are equivalent to
\be
\begin{array}{ll}
\Bigg. \langle A,P(a,d) \, b^m c^n \rangle = 
\displaystyle{\frac{\partial}{\partial a}} 
P(a,d) \Big\vert_{a=d=1} \ \delta_{m0} \delta_{n0} \,, &
\langle B,P(a,d) \, b^m c^n \rangle = P(a,d) \Big\vert_{a=d=1} 
\ \delta_{m1} \delta_{n0} \,, \cr
\Bigg. \langle C,P(a,d) \, b^m c^n \rangle = P(a,d) \Big\vert_{a=d=1} 
\ \delta_{m0} \delta_{n1} \,, &
\langle D,P(a,d) \, b^m c^n \rangle = 
\displaystyle{\frac{\partial}{\partial d}} 
P(a,d) \Big\vert_{a=d=1} \ \delta_{m0} \delta_{n0} \,. \cr
\end{array}
\ee
Therefore, the evaluation of the different (anti)commutators on $a^kd^l$
gives:
\begin{subequations}
\bea
\langle BC+CB,a^kd^l \rangle &=& \Delta^{kl}_{10,01}\tadp 
+ \Delta^{kl}_{01,10}\tadp \,, \label{prf2a} \\
\langle B^2,a^kd^l \rangle &=& \Delta^{kl}_{10,10}\tadp \,, \qquad
\langle C^2,a^kd^l \rangle = \Delta^{kl}_{01,01}\tadp \,, \label{prf2b} \\
\langle AB-BA,a^kd^l \rangle &=& \frac{\partial}{\partial a_1} 
\Delta^{kl}_{00,10}\tadp - \frac{\partial}{\partial a_2} 
\Delta^{kl}_{10,00}\tadp \,, \label{prf2c} \\
\langle AC-CA,a^kd^l \rangle &=& \frac{\partial}{\partial a_1} 
\Delta^{kl}_{00,01}\tadp - \frac{\partial}{\partial a_2} 
\Delta^{kl}_{01,00}\tadp \,, \\
\langle DB-BD,a^kd^l \rangle &=& \frac{\partial}{\partial d_1} 
\Delta^{kl}_{00,10}\tadp - \frac{\partial}{\partial d_2} 
\Delta^{kl}_{10,00}\tadp \,, \\
\langle DC-CD,a^kd^l \rangle &=& \frac{\partial}{\partial d_1} 
\Delta^{kl}_{00,01}\tadp - \frac{\partial}{\partial d_2} 
\Delta^{kl}_{01,00}\tadp \,, \label{prf2f} \\
\langle AD-DA,a^kd^l \rangle &=& \frac{\partial^2}{\partial a_1 \partial d_2} 
\Delta^{kl}_{00,00}\tadp - \frac{\partial^2}
{\partial a_2 \partial d_1} \Delta^{kl}_{00,00}\tadp \,. \label{prf2g}
\eea
\end{subequations}

\medskip

$\clubsuit$ 
We begin the proof by showing the following lemma:
\begin{lemap}\label{lemapA}
\beo
&& \Delta^{kl}_{10,01}\tadp = \frac14 \left( \frac{q^{2k+2l}-1}{q^2-1} 
+ \frac{q^{-2k-2l}-1}{q^{-2}-1} + 2(k-l) \right) \,, \\
&& \Delta^{kl}_{01,10}\tadp = \frac14 \left( \frac{q^{2k+2l}-1}{q^2-1} 
+ \frac{q^{-2k-2l}-1}{q^{-2}-1} - 2(k-l) \right) \,, \\
&& \Delta^{kl}_{01,01}\tadp = 
\Delta^{kl}_{10,10}\tadp = - \frac14 \left( \frac{q^{2k+2l}-1}{q^2-1} 
- \frac{q^{-2k-2l}-1}{q^{-2}-1} \right) \,. 
\eeo
\end{lemap}
Rewriting the multiplication law (\ref{multtwo}) in the following form:
\be
(b \pm c) \tad{$a$}{$d$} = M_\pm \tad{$a$}{$d$} (b \pm c) \quad\mbox{where}
\quad M_\pm = \left(\begin{array}{cc} r & \mp s \cr \pm s & -r \cr 
\end{array}\right) \,, \label{prf3}
\ee
it follows from eqs. (\ref{prf1}) and (\ref{prf3}) that
\bea
\Delta(a^{k+1}d^l) &=& (a \otimes a + b \otimes c) \Delta(a^kd^l) \nonumber \\
&=& (a \otimes a) \sum_{i,j,i',j' \in\{0,1\}} \Delta^{kl}_{ij,i'j'} \tadb 
b^ic^j \otimes b^{i'}c^{j'} \nonumber \\
&&+ \frac14 \ \sum_{i,j,i',j' \in\{0,1\}} \ \sum_{\eps_1,\eps_2=\pm 1} 
\eps_2 \Delta^{kl}_{ij,i'j'} \bigg(M_{\eps_1}\tad{$a_1$}{$d_1$},M_{\eps_2}
\tad{$a_2$}{$d_2$}\bigg) \nonumber \\
&& \hspace{60mm} (b^{i+1}c^j + \eps_1 b^ic^{j+1}) \otimes (b^{i'+1}c^{j'}
+ \eps_2 b^{i'}c^{j'+1}) \,, \label{prf25}
\eea
Looking at the different terms in $\Delta(a^{k+1}d^l)$, we get
\begin{subequations}
\label{prf4}
\bea
\Delta^{k+1,l}_{00,00}\tadb &=& a_1a_2 \ \Delta^{kl}_{00,00}\tadb 
+ \frac{s^2}{16} (a_1^2-d_1^2)(a_2^2-d_2^2) \sum_{\eps_1,\eps_2=\pm 1} 
(\eps_2 \Delta^{kl}_{10,10} + \Delta^{kl}_{10,01} \nonumber \\
&&+ \eps_1\eps_2 \Delta^{kl}_{01,10} + \eps_1 \Delta^{kl}_{01,01})
\bigg(M_{\eps_1}\tad{$a_1$}{$d_1$},M_{\eps_2}\tad{$a_2$}{$d_2$}\bigg) 
\,, \nonumber \\
&& \label{prf4a} \\
&& \nonumber \\
\Delta^{k+1,l}_{10,10}\tadb &=& a_1a_2 \ \Delta^{kl}_{10,10}\tadb 
+ \sum_{\eps_1,\eps_2=\pm 1} (\frac14 \eps_2 \Delta^{kl}_{00,00} 
+ \frac{s}{8} (a_1^2-d_1^2) \eps_1\eps_2 \Delta^{kl}_{11,00} \nonumber \\
&&+ \frac{s}{8} (a_2^2-d_2^2)\Delta^{kl}_{00,11} + \frac{s^2}{16} 
(a_1^2-d_1^2)(a_2^2-d_2^2) \eps_1 \Delta^{kl}_{11,11}) \bigg(M_{\eps_1}
\tad{$a_1$}{$d_1$},M_{\eps_2}\tad{$a_2$}{$d_2$}\bigg) \,, \nonumber \\
&& \\
&& \nonumber \\
\Delta^{k+1,l}_{01,01}\tadb &=& a_1a_2 \ \Delta^{kl}_{01,01}\tadb 
+ \sum_{\eps_1,\eps_2=\pm 1} (\frac14 \eps_1 \Delta^{kl}_{00,00} 
+ \frac{s}{8} (a_2^2-d_2^2) \eps_1\eps_2 \Delta^{kl}_{00,11} \nonumber \\
&&+ \frac{s}{8} (a_1^2-d_1^2)\Delta^{kl}_{11,00} + \frac{s^2}{16} 
(a_1^2-d_1^2)(a_2^2-d_2^2) \eps_2 \Delta^{kl}_{11,11}) \bigg(M_{\eps_1}
\tad{$a_1$}{$d_1$},M_{\eps_2}\tad{$a_2$}{$d_2$}\bigg) \,, \nonumber \\
&& \\
&& \nonumber \\
\Delta^{k+1,l}_{01,10}\tadb &=& a_1a_2 \ \Delta^{kl}_{01,10}\tadb 
+ \sum_{\eps_1,\eps_2=\pm 1} (\frac14 \eps_1\eps_2 \Delta^{kl}_{00,00} 
+ \frac{s}{8} (a_1^2-d_1^2) \eps_2 \Delta^{kl}_{11,00} \nonumber \\
&&+ \frac{s}{8} (a_2^2-d_2^2) \eps_1 \Delta^{kl}_{00,11} + \frac{s^2}{16} 
(a_1^2-d_1^2)(a_2^2-d_2^2) \Delta^{kl}_{11,11}) \bigg(M_{\eps_1}
\tad{$a_1$}{$d_1$},M_{\eps_2}\tad{$a_2$}{$d_2$}\bigg) \,, \nonumber \\
&& \\
&& \nonumber \\
\Delta^{k+1,l}_{10,01}\tadb &=& a_1a_2 \ \Delta^{kl}_{10,01}\tadb 
+ \sum_{\eps_1,\eps_2=\pm 1} (\frac14 \Delta^{kl}_{00,00} + \frac{s}{8} 
(a_1^2-d_1^2) \eps_1 \Delta^{kl}_{11,00} \nonumber \\
&&+ \frac{s}{8} (a_2^2-d_2^2) \eps_2 \Delta^{kl}_{00,11} + \frac{s^2}{16} 
(a_1^2-d_1^2)(a_2^2-d_2^2) \eps_1\eps_2 \Delta^{kl}_{11,11}) \bigg(M_{\eps_1}
\tad{$a_1$}{$d_1$},M_{\eps_2}\tad{$a_2$}{$d_2$}\bigg) \,. \nonumber \\
&& 
\eea
\end{subequations}
Taking now the values $a_1 = a_2 = d_1 = d_2 = 1$, one has easily
\begin{subequations}
\bea
&& \Delta^{k+1,l}_{10,01}\tadp = \Delta^{kl}_{10,01}\tadp + \frac14 
\sum_{\eps_1,\eps_2=\pm 1} \Delta^{kl}_{00,00}\tadi \,, \label{prf5a} \\
&& \Delta^{k+1,l}_{01,10}\tadp = \Delta^{kl}_{01,10}\tadp + \frac14 
\sum_{\eps_1,\eps_2=\pm 1} \eps_1\eps_2 \Delta^{kl}_{00,00}\tadi \,, \\
&& \Delta^{k+1,l}_{01,01}\tadp = \Delta^{kl}_{01,01}\tadp + \frac14 
\sum_{\eps_1,\eps_2=\pm 1} \eps_1 \Delta^{kl}_{00,00}\tadi \,, \\
&& \Delta^{k+1,l}_{10,10}\tadp = \Delta^{kl}_{10,10}\tadp + \frac14 
\sum_{\eps_1,\eps_2=\pm 1} \eps_2\Delta^{kl}_{00,00}\tadi \,, \label{prf5d} \\
&& \Delta^{k+1,l}_{00,00}\tadi = q^{-\eps_1-\eps_2} \Delta^{kl}_{00,00}\tadi 
\,. \label{prf5e}
\eea
\end{subequations}
The last equation (\ref{prf5e}) can be solved and one finds
\be
\Delta^{kl}_{00,00}\tadi = q^{-k(\eps_1+\eps_2)} 
\Delta^{0l}_{00,00}\tadi \,. \label{prf6}
\ee
In the same way, the recursion formulae for $l$ is given by
\bea
\Delta(a^kd^{l+1}) &=& (c \otimes b + d \otimes d) \Delta(a^kd^l) \nonumber \\
&=& (d \otimes d) \sum_{i,j,i',j' \in\{0,1\}} \Delta^{kl}_{ij,i'j'} \tadb 
b^ic^j \otimes b^{i'}c^{j'} \nonumber \\
&&+ \frac14 \ \sum_{i,j,i',j' \in\{0,1\}} \ \sum_{\eps_1,\eps_2=\pm 1} 
\eps_1 \Delta^{kl}_{ij,i'j'} \bigg(M_{\eps_1}\tad{$a_1$}{$d_1$},M_{\eps_2}
\tad{$a_2$}{$d_2$}\bigg) \nonumber \\
&& \hspace{60mm} (b^{i+1}c^j + \eps_1 b^ic^{j+1}) \otimes (b^{i'+1}c^{j'} 
+ \eps_2 b^{i'}c^{j'+1}) \,.
\eea
Looking at the different terms in $\Delta(a^kd^{l+1})$, we obtain relations
analogous to (\ref{prf4}), that lead to
\begin{subequations}
\bea
&& \Delta^{k,l+1}_{10,01}\tadp = \Delta^{kl}_{10,01}\tadp + \frac14 
\sum_{\eps_1,\eps_2=\pm 1} \eps_1\eps_2 \Delta^{kl}_{00,00}\tadi 
\,, \label{prf7a} \\
&& \Delta^{k,l+1}_{01,10}\tadp = \Delta^{kl}_{01,10}\tadp + \frac14 
\sum_{\eps_1,\eps_2=\pm 1} \Delta^{kl}_{00,00}\tadi \,, \\
&& \Delta^{k,l+1}_{01,01}\tadp = \Delta^{kl}_{01,01}\tadp + \frac14 
\sum_{\eps_1,\eps_2=\pm 1} \eps_2 \Delta^{kl}_{00,00}\tadi \,, \\
&& \Delta^{k,l+1}_{10,10}\tadp = \Delta^{kl}_{10,10}\tadp + \frac14 
\sum_{\eps_1,\eps_2=\pm 1} \eps_1 \Delta^{kl}_{00,00}\tadi 
\,, \label{prf7d} \\
&& \Delta^{k,l+1}_{00,00}\tadi = q^{-\eps_1-\eps_2} \Delta^{kl}_{00,00}\tadi 
\,. \label{prf7e}
\eea
\end{subequations}
Choosing $k=0$ in eq. (\ref{prf7e}) and taking into account eq. (\ref{prf6}), 
it follows that
\be
\Delta^{kl}_{00,00}\tadi = q^{-(k+l)(\eps_1+\eps_2)} \,.
\ee
Plugging this last result into eqs. (\ref{prf5a})--(\ref{prf5d}), 
one gets
\begin{subequations}
\bea
&& \Delta^{k+1,l}_{10,01}\tadp = \Delta^{kl}_{10,01}\tadp + \frac14 
(q^{2k+2l}+q^{-2k-2l}+2) \,, \\
&& \Delta^{k+1,l}_{01,10}\tadp = \Delta^{kl}_{01,10}\tadp + \frac14 
(q^{2k+2l}+q^{-2k-2l}-2) \,, \\
&& \Delta^{k+1,l}_{01,01}\tadp = \Delta^{kl}_{01,01}\tadp - \frac14 
(q^{2k+2l}-q^{-2k-2l}) \,, \\
&& \Delta^{k+1,l}_{10,10}\tadp = \Delta^{kl}_{10,10}\tadp - \frac14 
(q^{2k+2l}-q^{-2k-2l}) \,.
\eea
\end{subequations}
Hence using eqs. (\ref{prf7a})--(\ref{prf7d}), we obtain the results of 
the Lemma \ref{lemapA}. Then from eq. (\ref{prf2a}) it follows that
\be
\langle \{B,C\},a^kd^l \rangle = \frac12 \left(\frac{q^{2k+2l}-1}{q^2-1} + 
\frac{q^{-2k-2l}-1}{q^{-2}-1} \right) \,. \label{prf8}
\ee
Similarly eq. (\ref{prf2b}) leads to
\be
\langle B^2,a^kd^l \rangle = \langle C^2,a^kd^l \rangle = - \frac14 \left( 
\frac{q^{2k+2l}-1}{q^2-1} - \frac{q^{-2k-2l}-1}{q^{-2}-1} \right) 
\,. \label{prf9}
\ee

\medskip

$\clubsuit$ 
Along the same lines, one can derive recursion relations for the polynomials
$\Delta^{kl}_{ij,i'j'}$ where $i+j+i'+j'$ is odd -- this corresponds to the
choices $(ij,i'j') = (00,01)$, $(00,10), (01,00)$, $(10,00), (11,10)$, 
$(11,01), (10,11), (11,11)$. One has the two following lemmas:
\begin{lemap}\label{lemapB}
If $i+j+i'+j'$ is odd, the polynomials $\Delta^{kl}_{ij,i'j'}\tadp$ and
$\Delta^{kl}_{ij,i'j'}\tadi$ are identically vanishing.
\end{lemap}
\begin{lemap}\label{lemapC}
The derivatives 
$\displaystyle{\frac{\partial}{\partial a_1}\Delta^{kl}_{00,01}}$, 
$\displaystyle{\frac{\partial}{\partial a_1}\Delta^{kl}_{00,10}}$, 
$\displaystyle{\frac{\partial}{\partial a_2}\Delta^{kl}_{01,00}}$, 
$\displaystyle{\frac{\partial}{\partial a_2}\Delta^{kl}_{10,00}}$, 
$\displaystyle{\frac{\partial}{\partial d_1}\Delta^{kl}_{00,01}}$, 
$\displaystyle{\frac{\partial}{\partial d_1}\Delta^{kl}_{00,10}}$, 
$\displaystyle{\frac{\partial}{\partial d_2}\Delta^{kl}_{01,00}}$, 
$\displaystyle{\frac{\partial}{\partial d_2}\Delta^{kl}_{10,00}}$
taken at $\tad{$a_1$}{$d_1$},\tad{$a_2$}{$d_2$} = 
\tad{$1$}{$1$},\tad{$1$}{$1$}$ are all zero.
\end{lemap}
Consider for example the quantity $\Delta^{kl}_{00,01}$ which satisfies the
recursion relation
\bea
\Delta^{k+1,l}_{00,01}\tadb &=& a_1a_2 \ \Delta^{kl}_{00,01}\tadb 
+ (a_1^2-d_1^2) \sum_{\eps_1,\eps_2=\pm 1} (\frac{s}{8} \Delta^{kl}_{10,00} 
+ \frac{s}{8} \eps_1 \Delta^{kl}_{01,00} \nonumber \\
&& \hspace{-10mm} + \frac{s^2}{16} (a_2^2-d_2^2) \eps_2 \Delta^{kl}_{10,11} 
+ \frac{s^2}{16} (a_2^2-d_2^2) \eps_1 \eps_2 \Delta^{kl}_{01,11}) 
\bigg(M_{\eps_1}\tad{$a_1$}{$d_1$},M_{\eps_2}\tad{$a_2$}{$d_2$}\bigg) 
\,. \label{prf10}
\eea
obtained from eq. (\ref{prf25}). Thus for 
$\tad{$a_1$}{$d_1$},\tad{$a_2$}{$d_2$} = \tad{$1$}{$1$},\tad{$1$}{$1$}$, 
one has $\Delta^{k+1,l}_{00,01} = \Delta^{kl}_{00,01} = \Delta^{0l}_{00,01}$, 
while the recursion relation on $l$ leads to $\Delta^{0,l+1}_{00,01} = 
\Delta^{0l}_{00,01} = \Delta^{00}_{00,01}$; hence $\Delta^{kl}_{00,01}\tadp 
= 0$.
\\
For $\tad{$a_1$}{$d_1$},\tad{$a_2$}{$d_2$} = q^{-\eps_1}
\tad{\phantom{-}$1$}{-$1$},q^{-\eps_2}\tad{\phantom{-}$1$}{-$1$}$, one has
$\Delta^{k+1,l}_{00,01} = q^{-\eps_1-\eps_2} \Delta^{kl}_{00,01}$
and $\Delta^{0,l+1}_{00,01} = q^{-\eps_1-\eps_2} \Delta^{0l}_{00,01}$; 
hence $\Delta^{kl}_{00,01}\tadi = q^{-(k+l)(\eps_1+\eps_2)} 
\Delta^{00}_{00,01}\tadi = 0$.
\\
The same statement holds for the other cases, which proves Lemma \ref{lemapB}.
\\
Now taking the derivative with respect to $a_1$ of eq. (\ref{prf10}),
one gets
\bea
\frac{\partial}{\partial a_1} \Delta^{k+1,l}_{00,01}\tadp 
&=& \Delta^{kl}_{00,01} \tadp + \frac{\partial}
{\partial a_1} \Delta^{kl}_{00,01} \tadp \nonumber \\
&& + \frac{s}{4} \sum_{\eps_1,\eps_2=\pm 1} (\Delta^{kl}_{10,00} + \eps_1 
\Delta^{kl}_{01,00}) \tadi \,,
\eea
and analogous relations for the other combinations of the quadruplets
$(ij,i'j')$ with $i+j+i'+j'$ odd. From Lemma \ref{lemapB}, one has therefore
\be
\frac{\partial}{\partial a_1} \Delta^{kl}_{00,01} \tadp =
\frac{\partial}{\partial a_1} \Delta^{0l}_{00,01} \tadp \,.
\label{prf11}
\ee
Repeating the procedure for the recursion relations on $l$, one finds
that the r.h.s. of (\ref{prf11}) is zero and one concludes that 
\be
\frac{\partial}{\partial a_1} \Delta^{kl}_{00,01} \tadp = 0 \,.
\ee
The same statement holds for all the derivatives of the polynomials 
$\Delta^{kl}_{ij,i'j'}$ involved in eqs. (\ref{prf2c}) to (\ref{prf2f}),
which proves Lemma \ref{lemapC}. 
\\
It follows then from Lemmas \ref{lemapB} and \ref{lemapC} and eqs. 
(\ref{prf2c}) to (\ref{prf2f}) that
\be
\langle [A,B],a^kd^l \rangle = \langle [A,C],a^kd^l \rangle = 
\langle [D,B],a^kd^l \rangle = \langle [D,C],a^kd^l \rangle = 0 
\,. \label{prf12}
\ee
\medskip

$\clubsuit$ 
It remains to evaluate $\langle AD-DA,a^kd^l \rangle$. One has from eq.
(\ref{prf4a})
\bea
\frac{\partial^2}{\partial a_1 \partial d_2} \Delta^{k+1,l}_{00,00}\tadp 
&=& \frac{\partial^2}{\partial a_1 \partial d_2} \Delta^{kl}_{00,00}\tadp 
+ \frac{\partial}{\partial d_2} \Delta^{kl}_{00,00}\tadp \nonumber \\
&& \hspace{-30mm} - \frac{s^2}{4} \sum_{\eps_1,\eps_2 = \pm 1} 
(\Delta^{kl}_{10,01} + \eps_1 \Delta^{kl}_{01,01} + \eps_2 \Delta^{kl}_{10,10} 
+ \eps_1\eps_2 \Delta^{kl}_{01,10})\tadi \,, \\
\frac{\partial^2}{\partial a_2 \partial d_1} \Delta^{k+1,l}_{00,00}\tadp 
&=& \frac{\partial^2}{\partial a_2 \partial d_1} \Delta^{kl}_{00,00}\tadp 
+ \frac{\partial}{\partial d_1} \Delta^{kl}_{00,00}\tadp \nonumber \\
&& \hspace{-30mm} - \frac{s^2}{4} \sum_{\eps_1,\eps_2 = \pm 1} 
(\Delta^{kl}_{10,01} + \eps_1 \Delta^{kl}_{01,01} + \eps_2 \Delta^{kl}_{10,10} 
+ \eps_1\eps_2 \Delta^{kl}_{01,10})\tadi \,. 
\eea
Therefore
\bea
\left( \frac{\partial^2}{\partial a_1 \partial d_2} - \frac{\partial^2}
{\partial a_2 \partial d_1} \right) \Delta^{k+1,l}_{00,00}\tadp &=&
\left( \frac{\partial^2}{\partial a_1 \partial d_2} - \frac{\partial^2}
{\partial a_2 \partial d_1} \right) \Delta^{kl}_{00,00}\tadp \nonumber \\
&& + \left( \frac{\partial}{\partial d_2} - \frac{\partial}{\partial d_1} 
\right) \Delta^{kl}_{00,00}\tadp \,.
\eea
Then we use the following lemma:
\begin{lemap}\label{lemapD}
One has
\beo
&& \frac{\partial}{\partial a_1} \Delta^{kl}_{00,00}\tadp = 
\frac{\partial}{\partial a_2} \Delta^{kl}_{00,00}\tadp = k \,, \\
&& \frac{\partial}{\partial d_1} \Delta^{kl}_{00,00}\tadp = 
\frac{\partial}{\partial d_2} \Delta^{kl}_{00,00}\tadp = l \,.
\eeo
\end{lemap}
{From} equation (\ref{prf4a}), one has 
\be
\frac{\partial}{\partial d_i} \Delta^{k+1,l}_{00,00}\tadp = 
\frac{\partial}{\partial d_i} \Delta^{kl}_{00,00}\tadp 
\ee
and similarly
\be
\frac{\partial}{\partial d_i} \Delta^{0,l+1}_{00,00}\tadp = 
\Delta^{0l}_{00,00}\tadp + 
\frac{\partial}{\partial d_i} \Delta^{0l}_{00,00}\tadp \,.
\ee
Since $\Delta^{0,l+1}_{00,00}\tadp = \Delta^{0l}_{00,00}\tadp$
and $\Delta^{00}_{00,00}\tadp = 1$, one obtains the last line of Lemma
\ref{lemapD}. One gets the first line by exchanging the roles of $k$ and $l$.
\\
Then Lemma \ref{lemapD} implies
\be
\left( \frac{\partial^2}{\partial a_1 \partial d_2} - \frac{\partial^2}
{\partial a_2 \partial d_1} \right) \Delta^{k+1,l}_{00,00}\tadp =
\left( \frac{\partial^2}{\partial a_1 \partial d_2} - \frac{\partial^2}
{\partial a_2 \partial d_1} \right) \Delta^{kl}_{00,00}\tadp \,.
\ee
Similarly, one gets
\be
\left( \frac{\partial^2}{\partial a_1 \partial d_2} - \frac{\partial^2}
{\partial a_2 \partial d_1} \right) \Delta^{0,l+1}_{00,00}\tadp =
\left( \frac{\partial^2}{\partial a_1 \partial d_2} - \frac{\partial^2}
{\partial a_2 \partial d_1} \right) \Delta^{0l}_{00,00}\tadp \,.
\ee
Hence
\be
\left( \frac{\partial^2}{\partial a_1 \partial d_2} - \frac{\partial^2}
{\partial a_2 \partial d_1} \right) \Delta^{kl}_{00,00}\tadp =
\left( \frac{\partial^2}{\partial a_1 \partial d_2} - \frac{\partial^2}
{\partial a_2 \partial d_1} \right) \Delta^{00}_{00,00}\tadp = 0 \,. 
\ee
Therefore, from eq. (\ref{prf2g}), one obtains
\be
\langle [A,D],a^kd^l \rangle = 0 \,. \label{prf15}
\ee

\medskip

$\clubsuit$ 
Now we have to compute the evaluation of the (anti)commutators between
$A,B,C,D$ on the generic elements $a^kd^lb$ and $a^kd^lc$ of the
Poincar\'e--Birkhoff--Witt basis of $\cA$. Let us define
\be
\Delta(a^kd^lb) = \Delta(a^kd^l) (a \otimes b + b \otimes d) = 
\sum_{i,j,i',j' \in\{0,1\}} \beta^{kl}_{ij,i'j'} \tadb \ b^ic^j \otimes 
b^{i'}c^{j'} \,. \label{prf30}
\ee
It is not difficult to obtain the expressions of the polynomials
$\beta^{kl}_{ij,i'j'}$ in terms of the $\Delta^{kl}_{ij,i'j'}$'s from the
multiplication law (\ref{multtwo}).
One obtains (the other possibilities are of no interest for our goal):
\bea
\beta^{kl}_{00,00} &=& \shalf s a_1(a_2^2-d_2^2) \Delta^{kl}_{00,10} 
+ \shalf s d_2(a_1^2-d_1^2) \Delta^{kl}_{10,00} \,, \nonumber \\
\beta^{kl}_{10,00} &=& d_2 \Delta^{kl}_{00,00}
+ \shalf s ra_1(a_2^2-d_2^2) \Delta^{kl}_{10,10}
- \shalf s^2 d_1(a_2^2-d_2^2) \Delta^{kl}_{01,10} \,, \nonumber \\
\beta^{kl}_{01,00} &=& \shalf rs a_1(a_2^2-d_2^2) \Delta^{kl}_{01,10}
- \shalf s^2 d_1(a_2^2-d_2^2) \Delta^{kl}_{10,10} 
+ d_2(a_1^2-d_1^2) \Delta^{kl}_{11,00} \,, \nonumber \\
\beta^{kl}_{00,10} &=& a_1 \Delta^{kl}_{00,00} 
+ \shalf s^2 a_2(a_1^2-d_1^2) \Delta^{kl}_{10,01} 
- \shalf rs d_2(a_1^2-d_1^2) \Delta^{kl}_{10,10} \,, \nonumber \\
\beta^{kl}_{00,01} &=& \shalf s a_1(a_2^2-d_2^2)\Delta^{kl}_{00,11}
+ \shalf s^2 a_2(a_1^2-d_1^2)\Delta^{kl}_{10,10}
- \shalf rs d_2(a_1^2-d_1^2)\Delta^{kl}_{10,01} \,, \label{prf31} \\
\beta^{kl}_{10,10} &=& ra_1 \Delta^{kl}_{10,00} - sd_1 \Delta^{kl}_{01,00}
+ sa_2 \Delta^{kl}_{00,01} - rd_2 \Delta^{kl}_{00,10} \,, \nonumber \\
\beta^{kl}_{01,01} &=& \shalf rs a_1(a_2^2-d_2^2) \Delta^{kl}_{01,11}
- \shalf s^2 d_1(a_2^2-d_2^2) \Delta^{kl}_{10,11} 
+ \shalf s^2 a_2(a_1^2-d_1^2) \Delta^{kl}_{11,10}
- \shalf rs d_2(a_1^2-d_1^2) \Delta^{kl}_{11,01} \,, \nonumber \\
\beta^{kl}_{10,01} &=& sa_2 \Delta^{kl}_{00,10} - rd_2 \Delta^{kl}_{00,01}
+ \shalf rs a_1(a_2^2-d_2^2) \Delta^{kl}_{10,11}
- \shalf s^2 d_1(a_2^2-d_2^2) \Delta^{kl}_{01,11} \,, \nonumber \\
\beta^{kl}_{01,10} &=& ra_1 \Delta^{kl}_{01,00} - sd_1 \Delta^{kl}_{10,00}
+ \shalf s^2 a_2(a_1^2-d_1^2) \Delta^{kl}_{11,01}
- \shalf rs d_2(a_1^2-d_1^2) \Delta^{kl}_{11,10} \,. \nonumber
\eea
The evaluations of the (anti)commutators are given by
\begin{subequations}
\label{prf16}
\bea
\langle BC+CB,a^kd^lb \rangle &=& \beta^{kl}_{10,01} + \beta^{kl}_{01,10} 
\,, \\
\langle B^2,a^kd^lb \rangle &=& \beta^{kl}_{10,10} \,, \qquad
\langle C^2,a^kd^lb \rangle \ = \ \beta^{kl}_{01,01} \,, \\
\langle AB-BA,a^kd^lb \rangle &=& \frac{\partial}{\partial a_1} 
\beta^{kl}_{00,10} - \frac{\partial}{\partial a_2} \beta^{kl}_{10,00} 
\nonumber \\
&=& s^2(\Delta^{kl}_{10,01}+\Delta^{kl}_{01,10}) - 2rs\Delta^{kl}_{10,10} 
+ (1 + \frac{\partial}{\partial a_1} - \frac{\partial}{\partial a_2}) 
\Delta^{kl}_{00,00} \,, \\
\langle AC-CA,a^kd^lb \rangle &=& \frac{\partial}{\partial a_1} 
\beta^{kl}_{00,01} - \frac{\partial}{\partial a_2} \beta^{kl}_{01,00}
= 2s^2\Delta^{kl}_{10,10} - rs(\Delta^{kl}_{10,01}+\Delta^{kl}_{01,10}) \,, \\
\langle DB-BD,a^kd^lb \rangle &=& \frac{\partial}{\partial d_1} 
\beta^{kl}_{00,10} - \frac{\partial}{\partial d_2} \beta^{kl}_{10,00} 
\nonumber \\
&=& 2rs\Delta^{kl}_{10,10} - s^2(\Delta^{kl}_{10,01}+\Delta^{kl}_{01,10}) 
+ (-1 + \frac{\partial}{\partial d_1} - \frac{\partial}{\partial d_2})
\Delta^{kl}_{00,00} \,, \\
\langle DC-CD,a^kd^lb \rangle &=& \frac{\partial}{\partial d_1} 
\beta^{kl}_{00,01} - \frac{\partial}{\partial d_2} \beta^{kl}_{01,00} = 
rs(\Delta^{kl}_{10,01}+\Delta^{kl}_{01,10}) - 2s^2 \Delta^{kl}_{10,10} \,, \\
\langle AD-DA,a^kd^lb \rangle &=& \bigg(\frac{\partial^2}{\partial a_1 
\partial d_2} - \frac{\partial^2}{\partial a_2 \partial d_1} \bigg) 
\beta^{kl}_{00,00} \nonumber \\
&=& s(\Delta^{kl}_{10,00} - \Delta^{kl}_{00,10}
- (\frac{\partial}{\partial a_1} + \frac{\partial}{\partial d_1}) 
\Delta^{kl}_{10,00} 
+ (\frac{\partial}{\partial a_2} + \frac{\partial}{\partial d_2}) 
\Delta^{kl}_{00,10}) \,. 
\eea
\end{subequations}
where all polynomials $\beta^{kl}_{ij,i'j'}$, $\Delta^{kl}_{ij,i'j'}$ and 
their derivatives in (\ref{prf16}) are taken at 
$\tad{$a_1$}{$d_1$},\tad{$a_2$}{$d_2$} = \tad{$1$}{$1$},\tad{$1$}{$1$}$.

\medskip

Similarly, defining
\be
\Delta(a^kd^lc) = \Delta(a^kd^l) (c \otimes a + d \otimes c) = 
\sum_{i,j,i',j' \in\{0,1\}} \gamma^{kl}_{ij,i'j'} \tadb \ b^ic^j \otimes 
b^{i'}c^{j'} \,, \label{prf32}
\ee
one gets
\bea
\gamma^{kl}_{00,00} &=& \shalf s a_2(a_1^2-d_1^2) \Delta^{kl}_{01,00}
+ \shalf s d_1(a_2^2-d_2^2) \Delta^{kl}_{00,01} \,, \nonumber \\
\gamma^{kl}_{10,00} &=& \shalf s^2 a_1(a_2^2-d_2^2) \Delta^{kl}_{01,01}
- \shalf rs d_1(a_2^2-d_2^2) \Delta^{kl}_{10,01} 
+ \shalf s^2 a_2(a_1^2-d_1^2) \Delta^{kl}_{11,00} \,, \nonumber \\
\gamma^{kl}_{01,00} &=& a_2 \Delta^{kl}_{00,00} 
+ \shalf s^2 a_1(a_2^2-d_2^2) \Delta^{kl}_{10,01} 
- \shalf rs d_1(a_2^2-d_2^2) \Delta^{kl}_{01,01} \,, \nonumber \\
\gamma^{kl}_{00,10} &=& \shalf s^2 d_1(a_2^2-d_2^2) \Delta^{kl}_{00,11} 
+ \shalf rs a_2(a_1^2-d_1^2) \Delta^{kl}_{01,10} 
- \shalf s^2 d_2(a_1^2-d_1^2) \Delta^{kl}_{01,01} \,, \nonumber \\
\gamma^{kl}_{00,01} &=& d_1 \Delta^{kl}_{00,00} 
+ \shalf rs a_2(a_1^2-d_1^2) \Delta^{kl}_{01,01} 
- \shalf s^2 d_2(a_1^2-d_1^2) \Delta^{kl}_{01,10} \,, \label{prf33} \\
\gamma^{kl}_{10,10} &=& \shalf s^2 a_1(a_2^2-d_2^2) \Delta^{kl}_{01,11}
- \shalf rs d_1(a_2^2-d_2^2) \Delta^{kl}_{10,11} 
+ \shalf rs a_2(a_1^2-d_1^2) \Delta^{kl}_{11,10}
- \shalf s^2 d_2(a_1^2-d_1^2) \Delta^{kl}_{11,01} \,, \nonumber \\
\gamma^{kl}_{01,01} &=& a_1 \Delta^{kl}_{10,00} - rd_1 \Delta^{kl}_{01,00} 
+ ra_2 \Delta^{kl}_{00,01} - sd_2 \Delta^{kl}_{00,10} \,, \nonumber \\
\gamma^{kl}_{10,01} &=& \shalf s^2 a_1 \Delta^{kl}_{01,00}
- \shalf rs d_1 \Delta^{kl}_{10,00} 
+ \shalf rs a_2(a_1^2-d_1^2) \Delta^{kl}_{11,01}
- \shalf s^2 d_2(a_1^2-d_1^2) \Delta^{kl}_{11,10} \,, \nonumber \\
\gamma^{kl}_{01,10} &=& \shalf rs a_2 \Delta^{kl}_{00,10}
- \shalf s^2 d_2 \Delta^{kl}_{00,01} 
+ \shalf s^2 a_1(a_2^2-d_2^2) \Delta^{kl}_{10,11}
- \shalf rs d_1(a_2^2-d_2^2) \Delta^{kl}_{01,11} \,. \nonumber
\eea
Again, the evaluations of the (anti)commutators are given by
\begin{subequations}
\label{prf18}
\bea
\langle BC+CB,a^kd^lc \rangle &=& \gamma^{kl}_{10,01} + \gamma^{kl}_{01,10} 
\,, \\
\langle B^2,a^kd^lc \rangle &=& \gamma^{kl}_{10,10} \,, \qquad
\langle C^2,a^kd^lc \rangle \ = \ \gamma^{kl}_{01,01} \,, \\
\langle AB-BA,a^kd^lc \rangle &=& \frac{\partial}{\partial a_1} 
\gamma^{kl}_{00,10} - \frac{\partial}{\partial a_2} \gamma^{kl}_{10,00} = 
rs(\Delta^{kl}_{10,01}+\Delta^{kl}_{01,10}) - 2s^2\Delta^{kl}_{01,01} \,, \\
\langle AC-CA,a^kd^lc \rangle &=& \frac{\partial}{\partial a_1} 
\gamma^{kl}_{00,01} - \frac{\partial}{\partial a_2} \gamma^{kl}_{01,00} 
\nonumber \\
&=& 2rs\Delta^{kl}_{01,01} - s^2(\Delta^{kl}_{10,01} + \Delta^{kl}_{01,10}) 
+ (-1 +\frac{\partial}{\partial a_1} - \frac{\partial}{\partial a_2}) 
\Delta^{kl}_{00,00} \,, \\
\langle DB-BD,a^kd^lc \rangle &=& \frac{\partial}{\partial d_1} 
\gamma^{kl}_{00,10} - \frac{\partial}{\partial d_2} \gamma^{kl}_{10,00} = 
2s^2\Delta^{kl}_{01,01} - rs(\Delta^{kl}_{10,01}+\Delta^{kl}_{01,10}) \,, \\
\langle DC-CD,a^kd^lc \rangle &=& \frac{\partial}{\partial d_1} 
\gamma^{kl}_{00,01} - \frac{\partial}{\partial d_2} \gamma^{kl}_{01,00} 
\nonumber \\
&=& s^2(\Delta^{kl}_{10,01}+\Delta^{kl}_{01,10}) - 2rs\Delta^{kl}_{01,01}
+ (1 + \frac{\partial}{\partial d_1} - \frac{\partial}{\partial d_2}) 
\Delta^{kl}_{00,00} \,, \\
\langle AD-DA,a^kd^lc \rangle &=& \bigg(\frac{\partial^2}{\partial a_1 
\partial d_2} - \frac{\partial^2}{\partial a_2 \partial d_1} \bigg) 
\gamma^{kl}_{00,00} \nonumber \\
&=& s(\Delta^{kl}_{01,00} - \Delta^{kl}_{00,01}
- (\frac{\partial}{\partial a_1} + \frac{\partial}{\partial d_1}) 
\Delta^{kl}_{00,01}
+ (\frac{\partial}{\partial a_2} + \frac{\partial}{\partial d_2}) 
\Delta^{kl}_{01,00}) \,.
\eea
\end{subequations}
(the polynomials $\gamma^{kl}_{ij,i'j'}$, $\Delta^{kl}_{ij,i'j'}$ and their 
derivatives in (\ref{prf18}) are taken at 
$\tad{$a_1$}{$d_1$},\tad{$a_2$}{$d_2$} = \tad{$1$}{$1$},\tad{$1$}{$1$}$).

\medskip

Notice that all the polynomials $\Delta^{kl}_{ij,i'j'}$ and their derivatives
evaluated at $\tad{$a_1$}{$d_1$},\tad{$a_2$}{$d_2$} = 
\tad{$1$}{$1$},\tad{$1$}{$1$}$, that arise in the equations
(\ref{prf16}) and (\ref{prf18}), have
already been computed in the former steps of the proof (see in particular
Lemmas \ref{lemapA}, \ref{lemapB}, \ref{lemapC}, \ref{lemapD}).
It is therefore straightforward to obtain the evaluations of the 
(anti)commutators of $A,B,C,D$ on $a^kd^lb$ and $a^kd^lc$. One finds: 
\bea
&& \langle [A,C],a^kd^lb \rangle = \langle [D,B],a^kd^lc \rangle = 
- \squart (q^{2k+2l} - q^{-2k-2l}) \,, \nonumber \\
&& \langle [A,C],a^kd^lc \rangle = \langle [D,B],a^kd^lb \rangle = 
- \shalf - \squart (q^{2k+2l} + q^{-2k-2l}) \,, \nonumber \\
&& \langle [A,B],a^kd^lb \rangle = \langle [D,C],a^kd^lc \rangle = 
\shalf + \squart (q^{2k+2l} + q^{-2k-2l}) \,, \label{prf20} \\
&& \langle [A,B],a^kd^lc \rangle = \langle [D,C],a^kd^lb \rangle = 
\squart (q^{2k+2l} - q^{-2k-2l}) \,, \nonumber \\
&& \langle [A,D],a^kd^lb \rangle = \langle \{B,C\},a^kd^lb \rangle = 
\langle B^2,a^kd^lb \rangle = \langle C^2,a^kd^lb \rangle = 0 \,, \nonumber \\
&& \langle [A,D],a^kd^lc \rangle = \langle \{B,C\},a^kd^lc \rangle = 
\langle B^2,a^kd^lc \rangle = \langle C^2,a^kd^lc \rangle = 0 \,. \nonumber
\eea

\medskip

$\clubsuit$
Finally, one defines
\bea
&& \Delta(a^kd^lbc) = \Delta(a^kd^l) \Delta(bc) = \Delta(a^kd^l) 
(ad \otimes bc + bc \otimes ad + rac \otimes ab - rdb \otimes dc) \nonumber \\
&& \phantom{\Delta(a^kd^lbc)} = \sum_{i,j,i',j' \in\{0,1\}} \mu^{kl}_{ij,i'j'} 
\tadb \ b^ic^j \otimes b^{i'}c^{j'} \,.
\eea
One gets
\bea
\mu^{kl}_{00,00} &=& \squart s^2 a_1d_1(a_2^2-d_2^2)^2 \Delta^{kl}_{00,11}
+ \squart rs^2 a_1a_2(a_1^2-d_1^2)(a_2^2-d_2^2)\Delta^{kl}_{01,10} \nonumber \\
&& - \squart rs^2 d_1d_2(a_1^2-d_1^2)(a_2^2-d_2^2)\Delta^{kl}_{10,01}
+ \squart s^2 a_2d_2(a_1^2-d_1^2)^2\Delta^{kl}_{11,00} \,, \nonumber \\
\mu^{kl}_{10,00} &=& \shalf s a_2d_2(a_1^2-d_1^2)\Delta^{kl}_{01,00}
- \shalf rs^2 d_1a_2(a_2^2-d_2^2)\Delta^{kl}_{00,10} 
+ \shalf sr^2 d_1d_2(a_2^2-d_2^2)\Delta^{kl}_{00,01} \nonumber \\
&& + \squart rs^3 (a_1^2+d_1^2)(a_2^2-d_2^2)^2\Delta^{kl}_{01,11}
- \squart s^2(r^2+s^2) a_1d_1(a_2^2-d_2^2)^2\Delta^{kl}_{10,11} \nonumber \\
&& - \squart rs^3 a_1d_2(a_1^2-d_1^2)(a_2^2-d_2^2)\Delta^{kl}_{11,01} 
+ \squart r^2s^2 a_1a_2(a_1^2-d_1^2)(a_2^2-d_2^2)\Delta^{kl}_{11,10}
\,, \nonumber \\
\mu^{kl}_{01,00} &=& \shalf s a_2d_2(a_1^2-d_1^2)\Delta^{kl}_{10,00} 
+ \shalf sr^2 a_1a_2(a_2^2-d_2^2)\Delta^{kl}_{00,10} 
- \shalf rs^2 a_1d_2(a_2^2-d_2^2)\Delta^{kl}_{00,01} \nonumber \\
&& - \squart s^2(r^2+s^2) a_1d_1(a_2^2-d_2^2)^2\Delta^{kl}_{01,11}
+ \squart rs^3 (a_1^2+d_1^2)(a_2^2-d_2^2)^2\Delta^{kl}_{10,11} \nonumber \\
&& - \squart rs^3 d_1a_2(a_1^2-d_1^2)(a_2^2-d_2^2)\Delta^{kl}_{11,10}
+ \squart r^2s^2 d_1d_2(a_1^2-d_1^2)(a_2^2-d_2^2)\Delta^{kl}_{11,01}
\,, \nonumber \\
\mu^{kl}_{00,10} &=& - \shalf rs^2 d_1a_2(a_1^2-d_1^2)\Delta^{kl}_{10,00}
+ \shalf sr^2 a_1a_2(a_1^2-d_1^2)\Delta^{kl}_{01,00} 
+ \shalf s a_1d_1(a_2^2-d_2^2)\Delta^{kl}_{00,01} \nonumber \\
&& - \squart rs^3 a_1d_2(a_1^2-d_1^2)(a_2^2-d_2^2)\Delta^{kl}_{01,11}
+ \squart r^2s^2 d_1d_2(a_1^2-d_1^2)(a_2^2-d_2^2)\Delta^{kl}_{10,11} 
\nonumber \\
&& - \squart s^2(r^2+s^2) a_2d_2(a_1^2-d_1^2)^2\Delta^{kl}_{11,10}
+ \squart rs^3 (a_1^2-d_1^2)^2(a_2^2+d_2^2)\Delta^{kl}_{11,01} \,, \nonumber \\
\mu^{kl}_{00,01} &=& \shalf sr^2 d_1d_2(a_1^2-d_1^2)\Delta^{kl}_{10,00}
- \shalf rs^2 a_1d_2(a_1^2-d_1^2)\Delta^{kl}_{01,00}
+ \shalf s a_1d_1(a_2^2-d_2^2)\Delta^{kl}_{00,10} \nonumber \\
&& + \squart r^2s^2 a_1a_2(a_1^2-d_1^2)(a_2^2-d_2^2)\Delta^{kl}_{01,11}
- \squart rs^3 d_1a_2(a_1^2-d_1^2)(a_2^2-d_2^2)\Delta^{kl}_{10,11} \\
&& - \shalf s^2(r^2+s^2) a_2d_2(a_1^2-d_1^2)^2\Delta^{kl}_{11,01}
+ \squart rs^3 (a_1^2-d_1^2)^2(a_2^2+d_2^2)\Delta^{kl}_{11,10} \,, \nonumber \\
\mu^{kl}_{10,10} &=& \shalf rs a_1a_2(a_1^2-d_1^2)\Delta^{kl}_{11,00}
- \shalf rs d_1d_2(a_2^2-d_2^2)\Delta^{kl}_{00,11}
+ rs^2 (a_1^2a_2^2-d_1^2d_2^2)\Delta^{kl}_{01,01} \nonumber \\
&& - \shalf s(r^2+s^2) a_1d_1(a_2^2-d_2^2)\Delta^{kl}_{10,01}
- \shalf s(r^2+s^2) a_2d_2(a_1^2-d_1^2)\Delta^{kl}_{01,10} \,, \nonumber \\
\mu^{kl}_{01,01} &=& \shalf rs a_1a_2(a_2^2-d_2^2)\Delta^{kl}_{00,11}
- \shalf rs d_1d_2(a_1^2-d_1^2)\Delta^{kl}_{11,00}
+ rs (a_1^2a_2^2-d_1^2d_2^2)\Delta^{kl}_{10,10} \nonumber \\
&& - \shalf s(r^2+s^2) a_1d_1(a_2^2-d_2^2)\Delta^{kl}_{01,10} 
- \shalf s(r^2+s^2) a_2d_2(a_1^2-d_1^2)\Delta^{kl}_{10,01} \,, \nonumber \\
\mu^{kl}_{10,01} &=& - r d_1d_2\Delta^{kl}_{00,00}
+ \squart rs^2 a_1a_2(a_1^2-d_1^2)(a_2^2-d_2^2)\Delta^{kl}_{11,11}
+ rs^2 (a_1^2a_2^2-d_1^2d_2^2)\Delta^{kl}_{01,10} \nonumber \\
&& - \shalf s(r^2+s^2) a_1d_1(a_2^2-d_2^2)\Delta^{kl}_{10,10}
- \shalf s(r^2+s^2) a_2d_2(a_1^2-d_1^2)\Delta^{kl}_{01,01} \,, \nonumber \\
\mu^{kl}_{01,10} &=& r a_1a_2\Delta^{kl}_{00,00}
- \squart rs^2 d_1d_2(a_1^2-d_1^2)(a_2^2-d_2^2)\Delta^{kl}_{11,11}
+ rs^2 (a_1^2a_2^2-d_1^2d_2^2)\Delta^{kl}_{10,01} \nonumber \\
&& - \shalf s(r^2+s^2) a_2d_2(a_1^2-d_1^2)\Delta^{kl}_{10,10}
- \shalf s(r^2+s^2) a_1d_1(a_2^2-d_2^2)\Delta^{kl}_{01,01} \,. \nonumber
\eea
Once again, the evaluations of the (anti)commutators are given by
\begin{subequations}
\label{prf22}
\bea
\langle BC+CB,a^kd^lbc \rangle &=& \mu^{kl}_{10,01} + \mu^{kl}_{01,10} 
\,, \label{prf21a} \\
\langle B^2,a^kd^lbc \rangle &=& \mu^{kl}_{10,10} \,, \qquad
\langle C^2,a^kd^lbc \rangle \ = \ \mu^{kl}_{01,01} \,, \label{prf21b} \\
\langle AB-BA,a^kd^lbc \rangle &=& \frac{\partial}{\partial a_1} 
\mu^{kl}_{00,10} - \frac{\partial}{\partial a_2} \mu^{kl}_{10,00} \nonumber \\
&=& - rs^2 \Delta^{kl}_{10,00} + sr^2 \Delta^{kl}_{01,00} - rs^2 
\Delta^{kl}_{00,10} + sr^2 \Delta^{kl}_{00,01} \,, \label{prf21c} \\
\langle AC-CA,a^kd^lbc \rangle &=& \frac{\partial}{\partial a_1} 
\mu^{kl}_{00,01} - \frac{\partial}{\partial a_2} \mu^{kl}_{01,00} \nonumber \\
&=& sr^2 \Delta^{kl}_{10,00} - rs^2 \Delta^{kl}_{01,00} + sr^2 
\Delta^{kl}_{00,10} - rs^2 \Delta^{kl}_{00,01} \,, \\
\langle DB-BD,a^kd^lc \rangle &=& \frac{\partial}{\partial d_1} 
\mu^{kl}_{00,10} - \frac{\partial}{\partial d_2} \mu^{kl}_{10,00} \nonumber \\
&=& rs^2 \Delta^{kl}_{10,00} - sr^2 \Delta^{kl}_{01,00} + rs^2 
\Delta^{kl}_{00,10} - sr^2 \Delta^{kl}_{00,01} \,, \\
\langle DC-CD,a^kd^lbc \rangle &=& \frac{\partial}{\partial d_1} 
\mu^{kl}_{00,01} - \frac{\partial}{\partial d_2} \mu^{kl}_{01,00} \nonumber \\
&=& - sr^2 \Delta^{kl}_{10,00} + rs^2 \Delta^{kl}_{01,00} - sr^2 
\Delta^{kl}_{00,10} + rs^2 \Delta^{kl}_{00,01} \,, \label{prf21f} \\
\langle AD-DA,a^kd^lbc \rangle &=& \bigg(\frac{\partial^2}{\partial a_1 
\partial d_2} - \frac{\partial^2}{\partial a_2 \partial d_1} \bigg) 
\mu^{kl}_{00,00} \,. \label{prf21g}
\eea
\end{subequations}
(the polynomials $\mu^{kl}_{ij,i'j'}$, $\Delta^{kl}_{ij,i'j'}$ and their 
derivatives in (\ref{prf22}) are taken at 
$\tad{$a_1$}{$d_1$},\tad{$a_2$}{$d_2$} = \tad{$1$}{$1$},\tad{$1$}{$1$}$).
\\
Expressions (\ref{prf21a}) and (\ref{prf21b}) are obviously vanishing while 
expressions (\ref{prf21c}) to (\ref{prf21f}) are zero thanks to Lemma 
\ref{lemapB}. Finally, $\displaystyle{\frac{\partial^2}{\partial a_1 
\partial d_2}} \mu^{kl}_{00,00}\tadp = \displaystyle{\frac{\partial^2}
{\partial a_2 \partial d_1}} \mu^{kl}_{00,00}\tadp = 4r(\Delta^{kl}_{10,01} - 
\Delta^{kl}_{01,10})$, so that expression (\ref{prf21g}) is also vanishing. 
Therefore one has
\bea
&& \langle [A,B],a^kd^lbc\rangle = \langle [A,C],a^kd^lbc\rangle = 
\langle [D,B],a^kd^lbc\rangle = \langle [D,C],a^kd^lbc\rangle = 0 
\,, \nonumber \\
&& \langle [A,D],a^kd^lbc\rangle = \langle \{B,C\},a^kd^lbc\rangle = 
\langle B^2,a^kd^lbc\rangle = \langle C^2,a^kd^lbc\rangle = 0 \,. \label{prf21}
\eea

\medskip
$\clubsuit$
The last step of the proof consists to interpret the formulae (\ref{prf8}), 
(\ref{prf9}), (\ref{prf12}), (\ref{prf15}), (\ref{prf20}), (\ref{prf21}), 
that are evaluations of the (anti)commutators of the elements $A,B,C,D$ of 
$\cU$ onto generic elements of the Poincar\'e--Birkhoff--Witt basis of $\cA$, 
as abstract formulae defining the algebra given in Theorem 
\ref{thmalgtwo}. One has
\be
\langle (A+D)^n,a^kd^l \rangle = \langle \otimes_n (A+D),\Delta^{(n)}(a^kd^l) 
\rangle \,.
\ee
The generalization of the formula (\ref{prf1}) for the $n$-fold coproduct
reads as
\be
\Delta^{(n)}(a^kd^l) = \sum_{i_1,j_1,\dots,i_n,j_n\in\{0,1\}} 
{\Delta^{(n)}}^{kl}_{i_1j_1,\dots,i_nj_n}
\bigg(\tad{$a_1$}{$d_1$} \dots \tad{$a_n$}{$d_n$} \bigg)
\ b^{i_1}c^{j_1} \otimes \dots \otimes b^{i_n}c^{j_n} \,, \label{prf24}
\ee
where $a_i = 1 \otimes \dots \otimes 1 \otimes a \otimes 1 \dots \otimes 1$
and $a$ stands at the place $i$ of the tensor product, with a similar 
definition for $d_i$. Thus one has
\be
\langle \otimes_n (A+D),\Delta^{(n)}(a^kd^l) \rangle = 
\left\langle \otimes_n (A+D),{\Delta^{(n)}}^{kl}_{00,\dots,00}
\bigg(\tad{$a_1$}{$d_1$} \dots \tad{$a_n$}{$d_n$} \bigg) \right\rangle \,. 
\ee
Now the main observation is that the terms in (\ref{prf24}) coming from 
$b^2$ or $c^2$ cancel when evaluated on $A+D$ since $\langle A+D,a^kd^lb^2 
\rangle = \langle A+D,a^kd^lc^2 \rangle = \langle A+D,\shalf 
sa^kd^l(a^2-d^2) \rangle = 0$. It follows that the only 
relevant term of ${\Delta^{(n)}}^{kl}_{00,\dots,00}\bigg(\tad{$a_1$}{$d_1$} 
\dots \tad{$a_n$}{$d_n$}\bigg)$ is $a_1^kd_1^l \dots a_n^kd_n^l = a^kd^l 
\otimes \dots \otimes a^kd^l$. Therefore
\be
\langle (A+D)^n,a^kd^l \rangle = (k+l)^n \,,
\ee
{from} which we easily deduce
\be
\langle q^{A+D},a^kd^l \rangle = \langle K,a^kd^l \rangle = q^{k+l} 
\,. \label{prf23}
\ee
Moreover one has form eqs. (\ref{prf30}), (\ref{prf31}), (\ref{prf32}), 
(\ref{prf33}) and the previous results (note the shift in the exponential !):
\be
\langle q^{A+D-1}B,a^kd^lb \rangle = \langle q^{A+D-1}C,a^kd^lc \rangle 
= q^{k+l} \,. \label{prf26}
\ee
Then comparing eqs. (\ref{prf8}), (\ref{prf9}), (\ref{prf12}), (\ref{prf15}), 
(\ref{prf20}), (\ref{prf21}) with the formulae (\ref{prf23}) and 
(\ref{prf26}), Theorem \ref{thmalgtwo} immediately follows.
\hfill \rule{5pt}{5pt}

\end{document}